\documentclass{amsart}
\usepackage{amsthm,amssymb,amsmath,color}
\usepackage{pdfsync,epsfig}

\newcommand{\red}{\color{darkred}}

\newcommand{\clt}{central limit theorem}

\newcommand{\garch}{{\rm GARCH}$(1,1)$}

\newcommand{\sta}{St\u aric\u a}
\newcommand{\ex}{{\rm e}\,}

\newcommand{\asy}{asymptotic}

\newcommand{\ts}{time series}
\newcommand{\tsa}{\ts\ analysis}
\newcommand{\garchpq}{{\rm GARCH}$(p,q)$}
\definecolor{darkblue}{rgb}{.1, 0.1,.8}
\definecolor{darkgreen}{rgb}{0,0.8,0.2}
\definecolor{darkred}{rgb}{.8, .1,.1}

\newtheorem{lemma}{Lemma}[section]

\newtheorem{theorem}[lemma]{Theorem}
\newcommand{\Leb}{{\mathbb L}{\mathbb E}{\mathbb B}}

\newtheorem{proposition}[lemma]{Proposition}
\newtheorem{definition}[lemma]{Definition}
\newtheorem{corollary}[lemma]{Corollary}
\newtheorem{example}[lemma]{Example}
\newtheorem{exercise}[lemma]{Exercise}
\newtheorem{remark}[lemma]{Remark}
\newtheorem{fig}[lemma]{Figure}
\newtheorem{tab}[lemma]{Table}

\newcommand{\MC}{Markov chain}

\newcommand{\bth}{\begin{theorem}}
\newcommand{\ethe}{\end{theorem}}
\newcommand{\sv}{stochastic volatility}
\newcommand{\bre}{\begin{remark}\em }
\newcommand{\ere}{\end{remark}}
\newcommand{\fre}{frequenc}
\newcommand{\ble}{\begin{lemma}}
\newcommand{\ele}{\end{lemma}}
\newcommand{\sre}{stochastic recurrence equation}
\newcommand{\pp}{point process}
\newcommand{\bde}{\begin{definition}}
\newcommand{\ede}{\end{definition}}
\newcommand{\bco}{\begin{corollary}}
\newcommand{\eco}{\end{corollary}}

\newcommand{\bpr}{\begin{proposition}}
\newcommand{\epr}{\end{proposition}}

\newcommand{\bexer}{\begin{exercise}}
\newcommand{\eexer}{\end{exercise}}

\newcommand{\bexam}{\begin{example}}
\newcommand{\eexam}{\end{example}}

\newcommand{\bfi}{\begin{fig}}
\newcommand{\efi}{\end{fig}}

\newcommand{\btab}{\begin{tab}}
\newcommand{\etab}{\end{tab}}

\newcommand{\per}{periodogram}

\newcommand{\fidi}{finite-dimensional distribution}
\newcommand{\rv}{random variable}

\newcommand{\MDA}{{\rm MDA}}

\newcommand{\var}{{\rm var}}

\newcommand{\cov}{{\rm cov}}

\newcommand{\rhs}{right-hand side}
\newcommand{\df}{distribution function}

\newcommand{\beao}{\begin{eqnarray*}}
\newcommand{\eeao}{\end{eqnarray*}\noindent}

\newcommand{\beam}{\begin{eqnarray}}
\newcommand{\eeam}{\end{eqnarray}\noindent}

\newcommand{\beqq}{\begin{equation}}
\newcommand{\eeqq}{\end{equation}\noindent}

\newcommand{\bce}{\begin{center}}
\newcommand{\ece}{\end{center}}

\newcommand{\barr}{\begin{array}}
\newcommand{\earr}{\end{array}}

\newcommand{\stp}{\stackrel{P}{\rightarrow}}
\newcommand{\std}{\stackrel{d}{\rightarrow}}

\newcommand{\stv}{\stackrel{v}{\rightarrow}}
\newcommand{\stw}{\stackrel{w}{\rightarrow}}

\newcommand{\w}{\omega}

\newcommand{\nto}{n\to\infty}
\newcommand{\kto}{k\to\infty}
\newcommand{\xto}{x\to\infty}

\newcommand{\ov}{\overline}
\newcommand{\wt}{\widetilde}
\newcommand{\wh}{\widehat}
\newcommand{\vep}{\varepsilon}

\newcommand{\la}{\lambda}

\newcommand{\regvary}{regularly varying}
\newcommand{\slvary}{slowly varying}
\newcommand{\regvar}{regular variation}

\newcommand{\bbr}{{\mathbb R}}

\newcommand{\bbz}{{\mathbb Z}}

\newcommand{\bbs}{{\mathbb S}}

\newcommand{\BM}{Brownian motion}

\newcommand{\con}{convergence}

\newcommand{\evt}{extreme value theory}
\newcommand{\evd}{extreme value distribution}

\newcommand{\st}{such that}
\newcommand{\fif}{if and only if}
\newcommand{\wrt}{with respect to}

\newcommand{\fct}{function}

\newcommand{\ds}{distribution}

\newcommand{\rep}{representation}
\newcommand{\cmt}{continuous mapping theorem}
\newcommand{\seq}{sequence}

\newcommand{\ins}{insurance}

\newcommand{\pro}{probabilit}

\newcommand{\ms}{measure}

\newcommand{\ld}{large deviation}
\newcommand{\bfx}{{\bf x}}

\newcommand{\bfy}{{\bf y}}

\newcommand{\bft}{{\bf t}}

\newcommand{\bfs}{{\bf s}}

\newcommand{\acf}{autocorrelation \fct }

\begin{document}
\today
\title[Measures of serial extremal dependence and their estimation]
{Measures of serial extremal dependence and their estimation}
\thanks{Thomas Mikosch's research is partly supported
by the Danish Research Council (FNU) Grants
09-072331 "Point process modelling and statistical inference"
and 10-084172 ``Heavy tail phenomena: Modeling and estimation''.
The research of Yuwei Zhao is supported by the Danish Research Council
Grant 10-084172. Richard A. Davis's research is supported by 
a Villum Kann Ramussen Research Grant 2011-2012, 
and NSF Grant DMS-1107031.}
\author[R.A Davis]{Richard A. Davis}
\author[T. Mikosch]{Thomas Mikosch}
\author[Y. Zhao]{Yuwei Zhao}
\address{Richard A. Davis, Department of Statistics,
Columbia University,
1255 Amsterdam Ave.
New York, NY 10027, U.S.A.}
\email{rdavis@stat.columbia.edu}

\address{Thomas Mikosch, University of Copenhagen, Department of Mathematics,
Universitetsparken 5,
DK-2100 Copenhagen\\ Denmark} \email{mikosch@math.ku.dk}
\address{Yuwei Zhao,
  University of Copenhagen, Department of Mathematics,
Universitetsparken 5,
DK-2100 Copenhagen\\ Denmark}           
            \email{zywmar@gmail.com}
\begin{abstract}
The goal of this paper is two-fold: 1. We review classical 
and recent \ms s
  of serial extremal dependence in a strictly stationary \ts\ as well
as their estimation. 2. We discuss recent concepts of heavy-tailed
\ts , including \regvar\ and max-stable processes.

Serial extremal dependence is typically characterized by clusters of
exceedances of high thresholds in the series. We start by discussing  
the notion of
extremal index of a univariate \seq , i.e. 
the reciprocal of the expected cluster
size, which has attracted major attention in the extremal value
literature.  Then we continue by introducing the extremogram which is
an \asy\ \acf\ for \seq s of extremal events in a \ts . In this
context, we discuss \regvar\ of a  \ts . This notion has been
useful for describing serial extremal dependence and heavy tails 
in a strictly stationary \seq . We briefly discuss the tail process
coined by Basrak and Segers to describe the dependence 
structure of \regvary\ \seq s in a probabilistic way. Max-stable
processes with Fr\'echet marginals are an important class of \regvary\
\seq s. Recently, this class has attracted attention 
for modeling and statistical purposes. We apply the extremogram to
max-stable processes. Finally, we discuss estimation of
the extremogram both in the time and  \fre y domains. 
\end{abstract}
\maketitle

\section{Introduction}\label{sec:introd}
\setcounter{equation}{0}
Measuring and estimating extremal dependence in a \ts\ is a rather challenging
problem. Since many real-life \ts , especially those arising in
finance and environmental applications,  are non-Gaussian
their dependence structure is not determined by their autocorrelation \fct .
Correlations are moments of the observations and as such
not well suited for describing the dependence of extremes which
typically arise from the tails of the underlying \ds .  
\subsection{The extremal index as reciprocal of the expected extremal
  cluster size}
Extremal dependence in a real-valued strictly stationary \seq\ $(X_t)$
can be described by the phenomenon of extremal clustering. Given some
sufficiently high threshold $u=u_n$, we would expect that exceedances
of this threshold should occur according to a homogeneous Poisson process. if $(X_t)$ is iid. On the other hand,
for dependent $(X_t)$ exhibiting extremal dependence, 
exceedances of $u$  should cluster in the sense that 
an exceedance of a high threshold is likely to be surrounded by neighboring
observations that also exceed the threshold.
Although the notion of extremal clustering is intuitively
appealing, a precise formulation is not so easy.
\par
The intuition about extremal clusters in a \ts\ 
can be made precise using \pp\ theory. In the classical monograph by
Leadbetter, Lindgren and Rootz\'en \cite{leadbetter:lindgren:rootzen:1983}
the {\em \pp\ of exceedances} of $u$ was used to describe clusters of 
extremes as an \asy\ phenomenon when the threshold $u_n$ converges to
the right endpoint of the \ds\ $F$ of $X$. (Here and in what follows,
$Y$ denotes a generic element of any strictly stationary \seq\
$(Y_t)$.) To be more precise, $(u_n)$ has to satisfy the condition 
$n\,\ov F(u_n)=n\,(1-F(u_n))\to \tau$ for some $\tau\in (0,\infty)$. Under this condition
and  mixing assumptions, the \pp es of exceedances converge weakly to
a compound Poisson process (see Hsing et
al. \cite{hsing:husler:leadbetter:1988}): 
\beam\label{eq:aa}
N_n=\sum_{i=1}^n \vep_{i/n} \,I_{\{X_i>u_n\}}\std N=\sum_{i=1}^\infty 
\xi _i \,\vep_{\Gamma_i}\,,
\eeam 
where the state space of the \pp es is $(0,1]$,
the points $0<\Gamma_1<\Gamma_2<\cdots$ constitute a homogeneous
Poisson process with intensity $\theta \tau$ on $(0,1]$ which is 
independent of an iid positive integer-valued \seq\ $(\xi_i)$. Here $\theta\in
[0,1]$ is the {\em extremal index} of the \seq\ $(X_t)$.
Thus, in an \asy\ way,
a cluster of extremes is located at the Poisson points $\Gamma_i$ with
corresponding size $\xi_i$. The cluster size \ds\ $P(\xi=k)$,
$k\ge 1$, contains plenty of information about the \ds\ of the
extremal clusters. However, most attention has been given to determine
the expected cluster size $E\xi$ which can be interpreted as 
reciprocal of $\theta$ as the following heuristic argument illustrates.
Applying \eqref{eq:aa} on the set $(0,1]$ and taking 
expectations on both sides of the limit relation, we observe that
\beao
EN_n(0,1]= n\ov F(u_n)\to \tau = EN(0,1]= 
E\xi \,E\# \{i\ge 1: \Gamma_i\le 1\}= E\xi \,(\theta\, \tau)\,, \quad
\nto .
\eeao
Thus $\theta= 1/E\xi$ with the convention that $E\xi=\infty$ for
$\theta=0$.
In the case of an iid \seq ,
$\xi_i\equiv 1$ a.s., i.e. $N$ collapses to a homogeneous
Poisson process and  $\theta=1$.
\par
Writing $M_n=\max(X_1,\ldots,X_n)$, we also observe that
\beao
P(M_n\le u_n)=P(N_n(0,1]=0)\to P(N(0,1]=0)=P(\#\{i\ge 1: \Gamma_i\le
1\}=0)
=\ex^{-\theta\,\tau}\,,
\eeao while for an iid \seq\ $(\wt X_t)$ with the same marginal \ds\ $F$
as for $(X_t)$ and $\wt M_n=\max(\wt X_1,\ldots,\wt X_n)$ we have 
\beao P(\wt M_n\le
u_n)= F^n(u_n)=\ex^{-n\ov F(u_n)(1+o(1))}\to \ex^{-\tau}\,.
\eeao
If $F$ belongs to the {\em maximum domain of attraction
of an extreme value \ds } $H$ ($F\in\MDA(H)$) there exist constants $c_n>0,
d_n\in\bbr$, $n\ge 1$, \st\ $P(c_n^{-1} (\wt M_n-d_n)\le x)\to H(x)$
for every $x\in {\rm supp}( H)$ (the support of $H$); cf.~Embrechts et
al. \cite{embrechts:kluppelberg:mikosch:1997}, Chapter 3.
Thus, 
writing $u_n(x)=c_n \,x + d_n$
and $\tau=\tau(x)=-\log H(x)$ for $x\in {\rm supp } (H)$, the
existence of an extremal value index $\theta$ of the \seq\ $(X_t)$
implies that
\beam\label{eq:conv}
P(c_n^{-1}(M_n- d_n)\le x)\to H^\theta(x)\,,\quad x\in\bbr\,.  
\eeam
\par
The concrete form of the extremal index is known for various standard
\ts\ models, including linear processes with iid 
subexponential noise (cf. \cite{embrechts:kluppelberg:mikosch:1997}, 
Section~5.5), Markov processes (see Leadbetter and Rootz\'en 
\cite{leadbetter:rootzen:1988}, Perfekt \cite{perfekt:1994})
and financial \ts\ models  such as GARCH 
(generalized autoregressive conditionally heteroscedastic) and SV
(stochastic volatility) models; cf. 
\cite{davis:mikosch:2001,davis:mikosch:2009a,davis:mikosch:2009b}.
Expressions of the extremal index for \regvary\ \seq s $(X_t)$
(see Section~\ref{subsec:regvar} for a definition)
in terms of the  points of the limiting \pp\ were given in 
Davis and Hsing \cite{davis:hsing:1995} and in terms of the limiting 
tail process in Basrak and Segers 
\cite{basrak:segers:2009}; see \eqref{eq:basrak} below. However, for most models 
these concrete expressions of $\theta$ are too complex 
to be useful in practice.
\par
An exception are Gaussian stationary \seq s $(X_t)$. Writing 
$\gamma_X(h)=\cov(X_0,X_h)$, $h\ge 0$, for the covariance \fct\ of
$(X_t)$, this \seq\ has extremal index $\theta=1$ under the very weak
condition $\gamma_X(h)= o(1/\log h)$ as $h\to\infty$ 
(so-called {\em Berman's condition}); see Leadbetter et al. 
\cite{leadbetter:lindgren:rootzen:1983}, cf. Theorem 4.4.8 in 
Embrechts et al. \cite{embrechts:kluppelberg:mikosch:1997}.
Notice that Berman's condition is satisfied for fractional Gaussian
noise and fractional Gaussian ARIMA processes (see Chapter 7 in Samorodnitsky
and Taqqu \cite{samorodnitsky:taqqu:1994}, and Section 13.2 
in Brockwell and Davis 
\cite{brockwell:davis:1991}). Subclasses of the latter processes
exhibit long range dependence  in the sense 
that $\sum_h|\gamma_X(h)|=\infty$.\footnote{This remark also indicates that 
long range dependence for extremes should not be defined via the
covariance \fct\ $\gamma_X$.  As explained above (see \eqref{eq:conv})
the existence of a positive extremal index $\theta$ ensures that the
type of the limiting \evd\ $H$ remains the same as in the iid
case. This is easily checked since the only possible
non-degenerate limit \ds s $H$ are the types of the {\em Fr\' echet \ds } 
$\Phi_\alpha(x)=\ex^{-x^{-\alpha}}$, $x,\alpha>0$, the {\em Weibull
  \ds } 
$\Psi_\alpha(x)=\ex^{-(-x)^\alpha}$, $x<0,\alpha>0$, and the {\em Gumbel
  \ds } $\Lambda(x)=\ex^{-\ex^{-x}}$, $x\in\bbr$. This is a con\seq\
of the Fisher-Tippett theorem; cf. 
Embrechts et al. \cite{embrechts:kluppelberg:mikosch:1997}, Theorem
3.2.3. The notion of long range dependence in an extreme value sense
would be reasonable if in \eqref{eq:conv} a limit \ds\ occurred which 
does not belong to the type of any of the 
three mentioned standard \evd s. This, 
however, can only be expected if a given stationary \seq\ $(X_t)$ with
$F\in\MDA(H)$ does not have an extremal index or if $\theta=0$. 
Examples of \seq s with zero extremal index are given  in 
Leadbetter et al.
\cite{leadbetter:lindgren:rootzen:1983} and 
Leadbetter \cite{leadbetter:1983}, but such examples are often
considered pathological; see also the discussion in 
Samorodnitsky \cite{samorodnitsky:2004} who
studied  infinite variance stable stationary \seq s with zero extremal
index and the boundary between short and long range extremal dependence for
these \seq .} 
\par
We conclude that any Gaussian stationary \seq s
which are
relevant for applications do not exhibit extremal clustering in the
sense that $\theta=1$. If $\theta=1$ one often says that
$(X_t)$ exhibits \asy\
independence of its extremes. However, the notion of 
{\em \asy\ independence } is not well defined and 
may have rather different meanings in the extreme value context, 
as we will observe later.
\par
Due to the complexity of expressions for the extremal index
it has been recognized early on that $\theta$ needs to be estimated
from real-life or simulated data. Various estimators were proposed in
the literature. Among them, the blocks and runs estimators are 
the most popular ones. These estimators are non-parametric estimators
of $\theta$ which, in different ways, define and count clusters in the 
sample and use this information to build estimators of $\theta$
under mixing conditions. In addition to the delicate choice of a
threshold $u_n$, these estimation techniques also involve the
construction of blocks of constant (but increasing with the sample
size $n$) length or of flexible length
depending on the local extremal behavior. These estimators
often exhibit a rather large uncertainty. 
\par
\begin{figure}[htbp]
\centerline{
\epsfig{figure=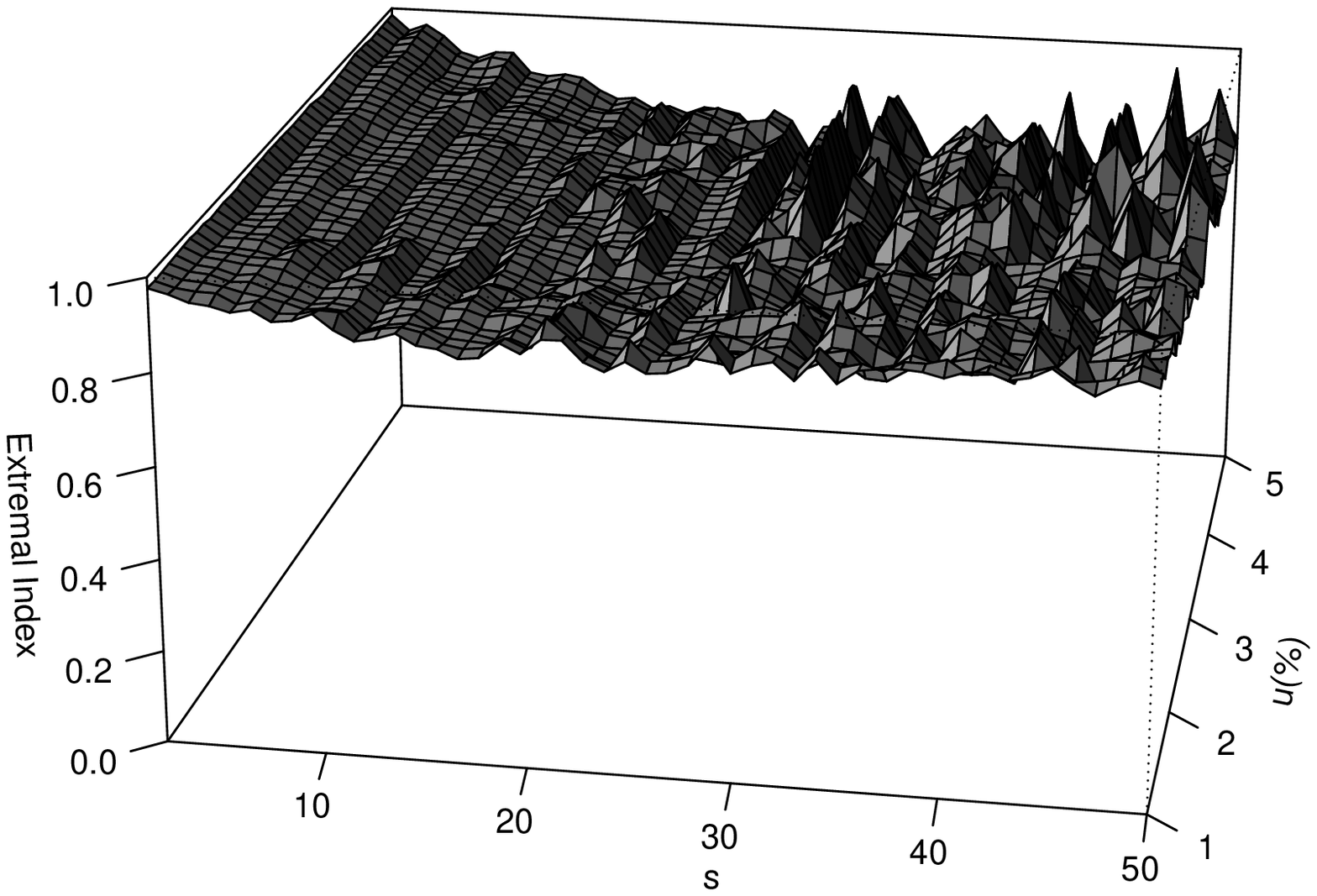,height=8cm,width=10cm}}
\centerline{
\epsfig{figure=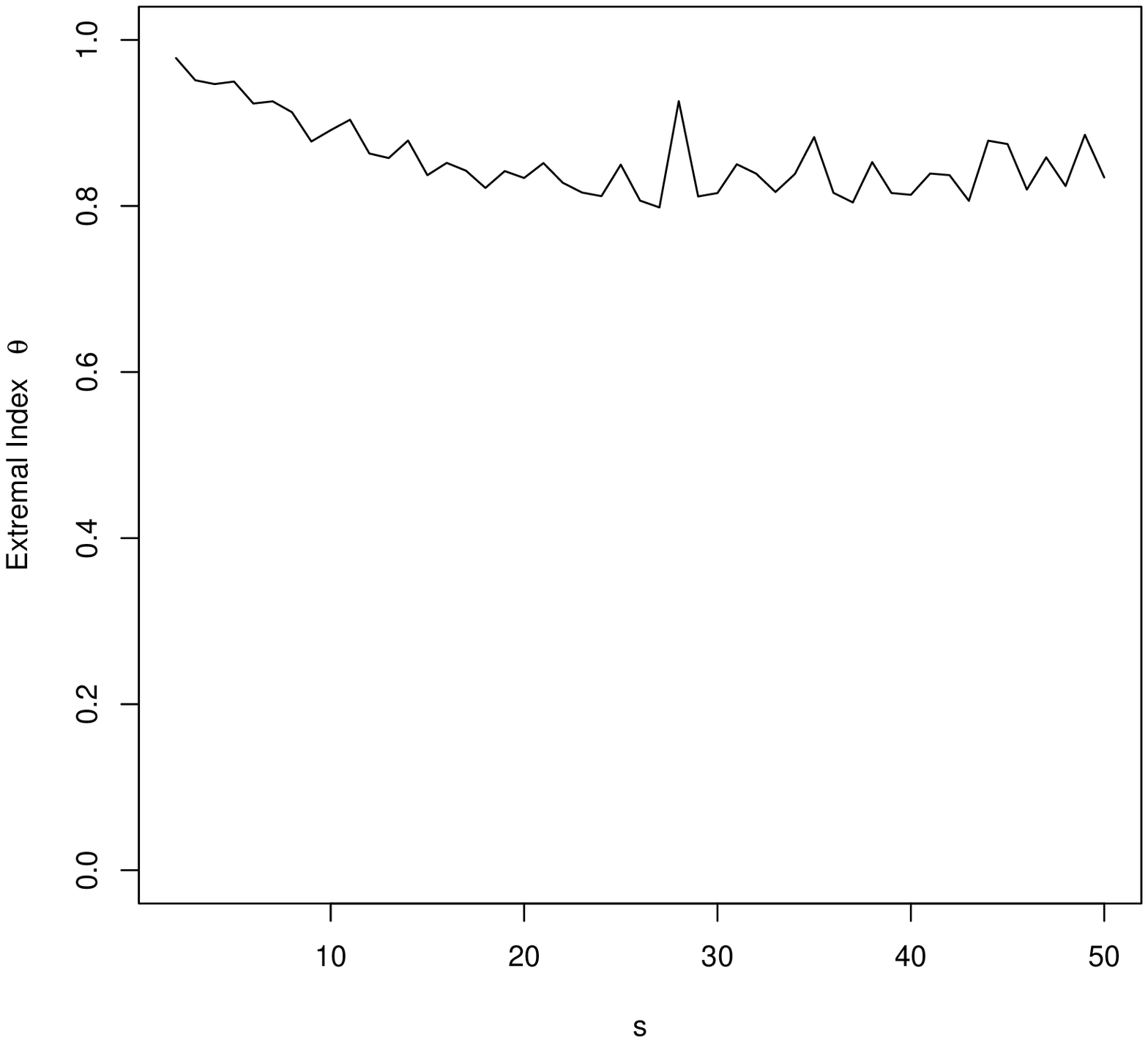,height=6cm,width=6cm,angle=0}
\epsfig{figure=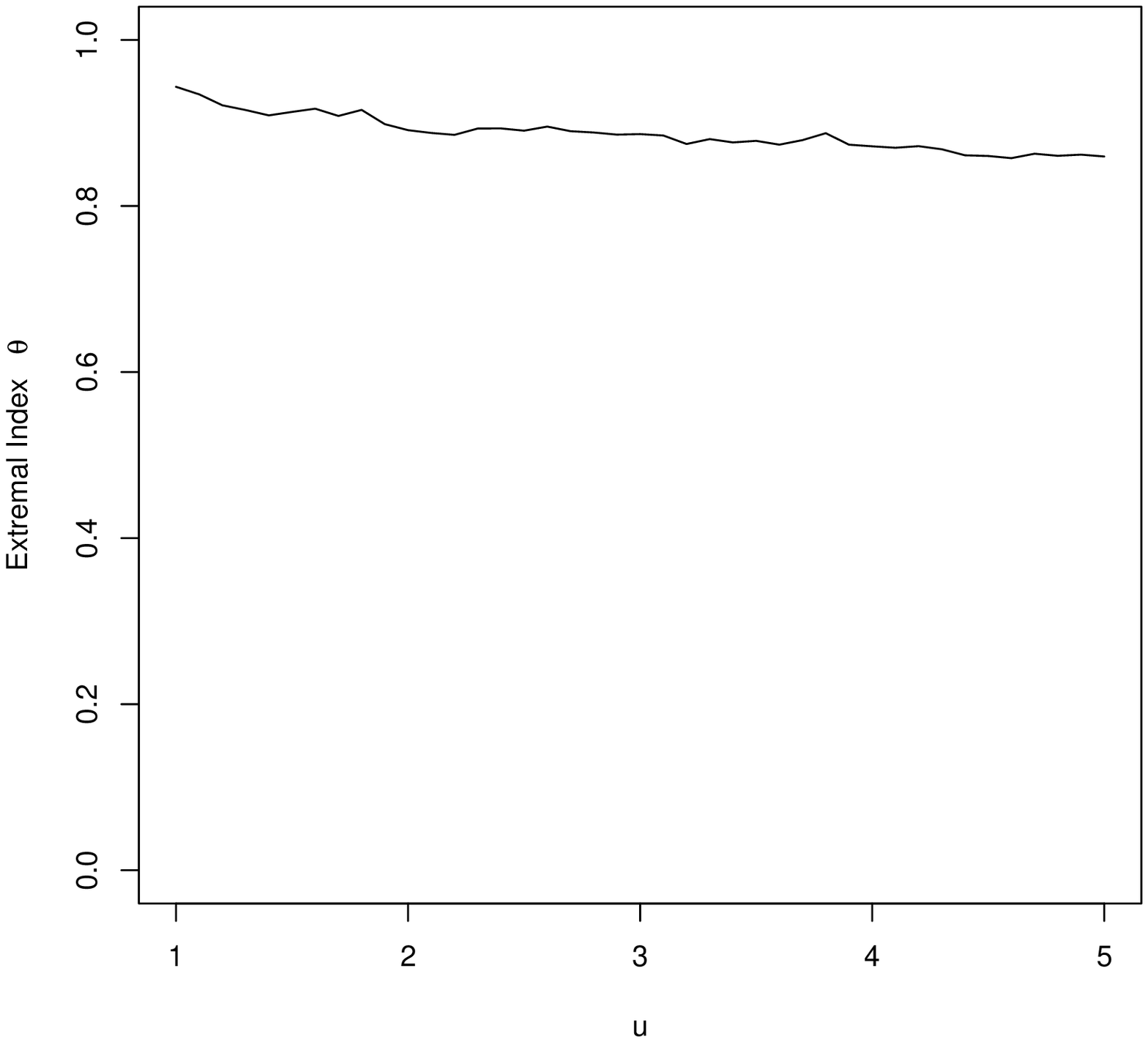,height=6cm,width=6cm,angle=0}
}
\bfi{\em Blocks estimator $\wh
  \theta$ of the extremal index $\theta$ for $31,757$ $5$-minute
  log-returns of Bank of America stock prices. The blocks
  estimator as a \fct\ of the block size $s$ and the $u\%$ upper order statistics
(top), for fixed $u=1.9\%$ and running $s$ (bottom left) and for fixed
$s=10$ and running $u$ (bottom right).  
 }\label{fig:index1}
\efi
\end{figure}
\begin{figure}[htbp]
\centerline{
\epsfig{figure=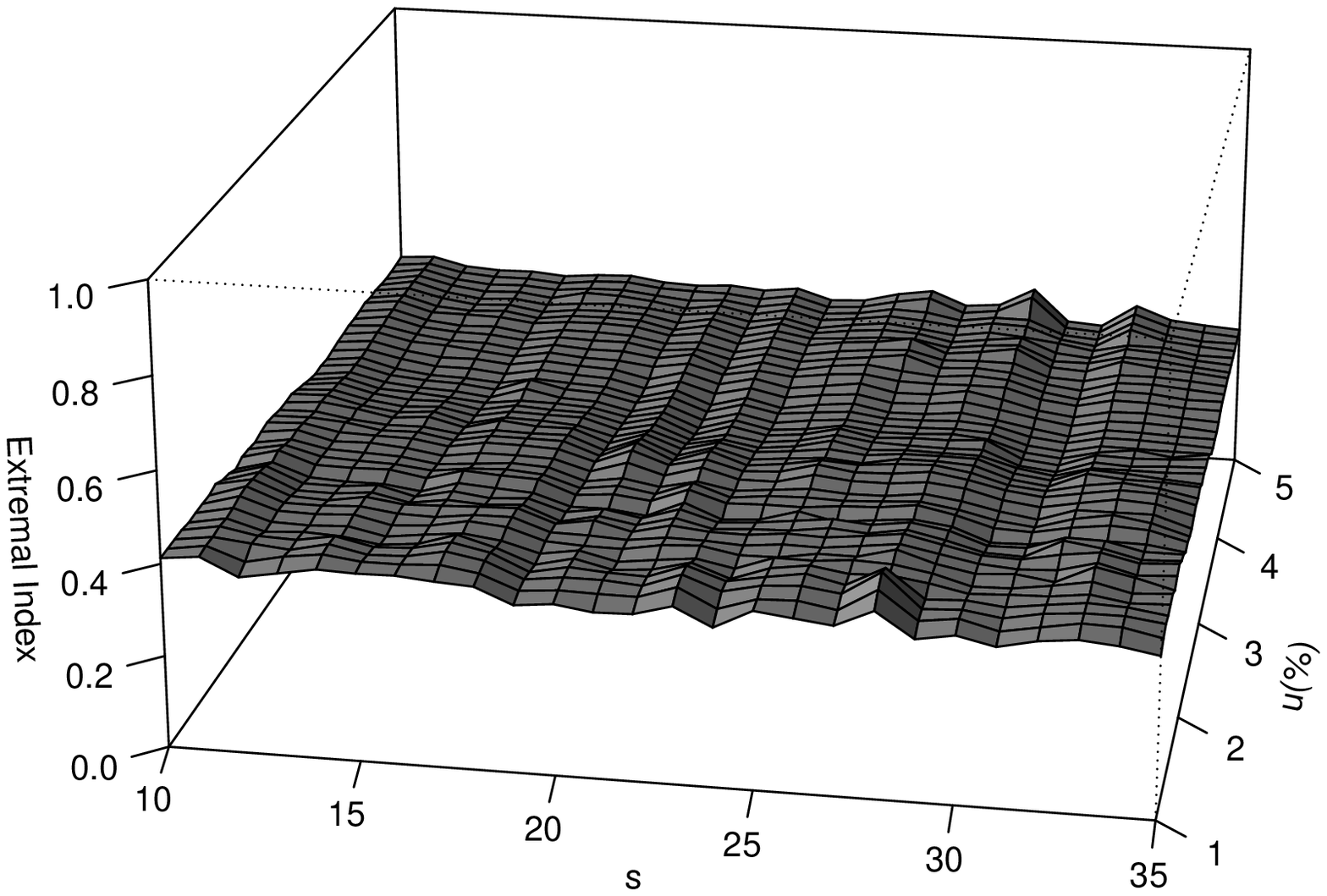,height=8cm,width=10cm,angle=0}}
\centerline{
\epsfig{figure=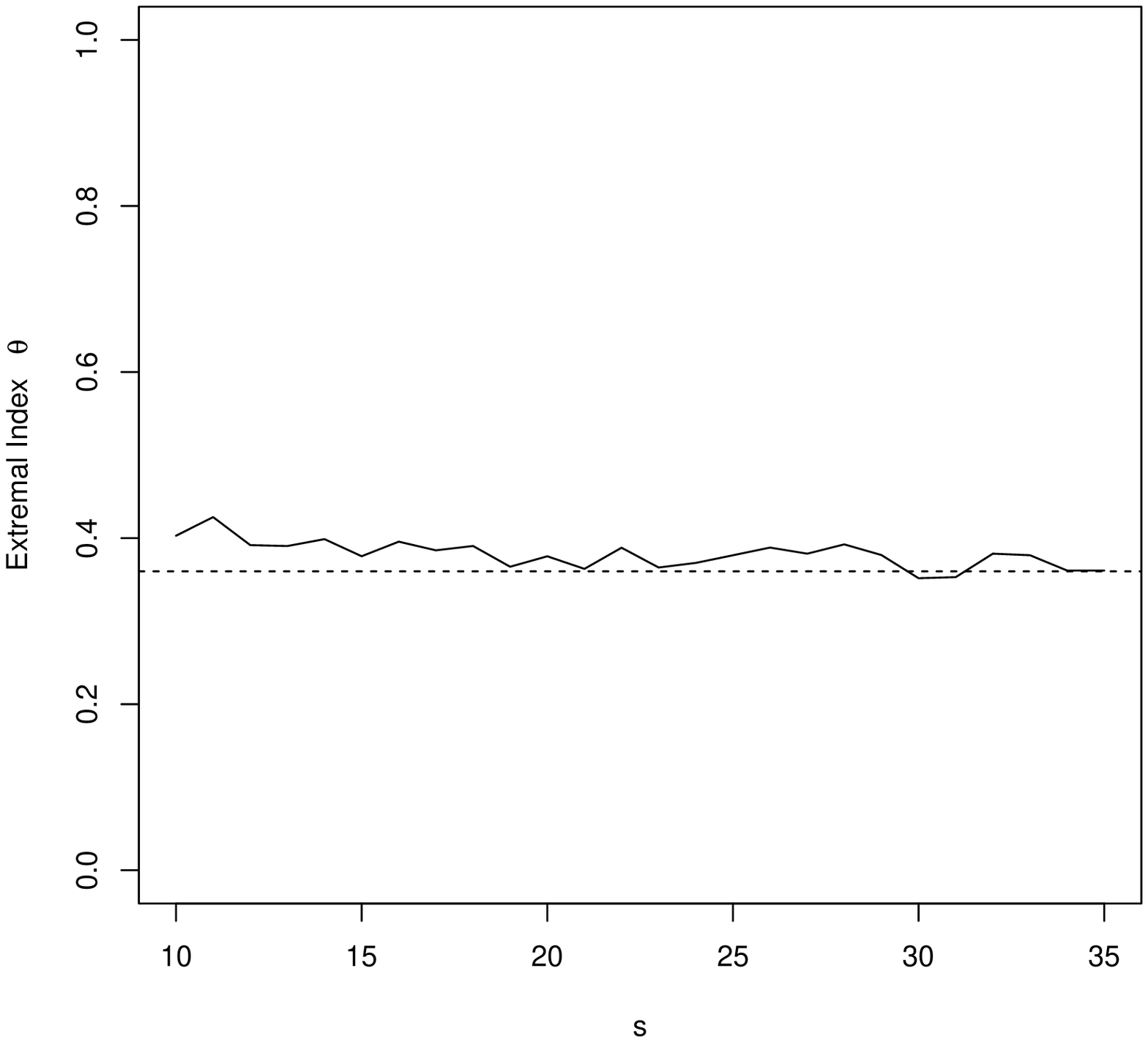,height=6cm,width=6cm,angle=0}
\epsfig{figure=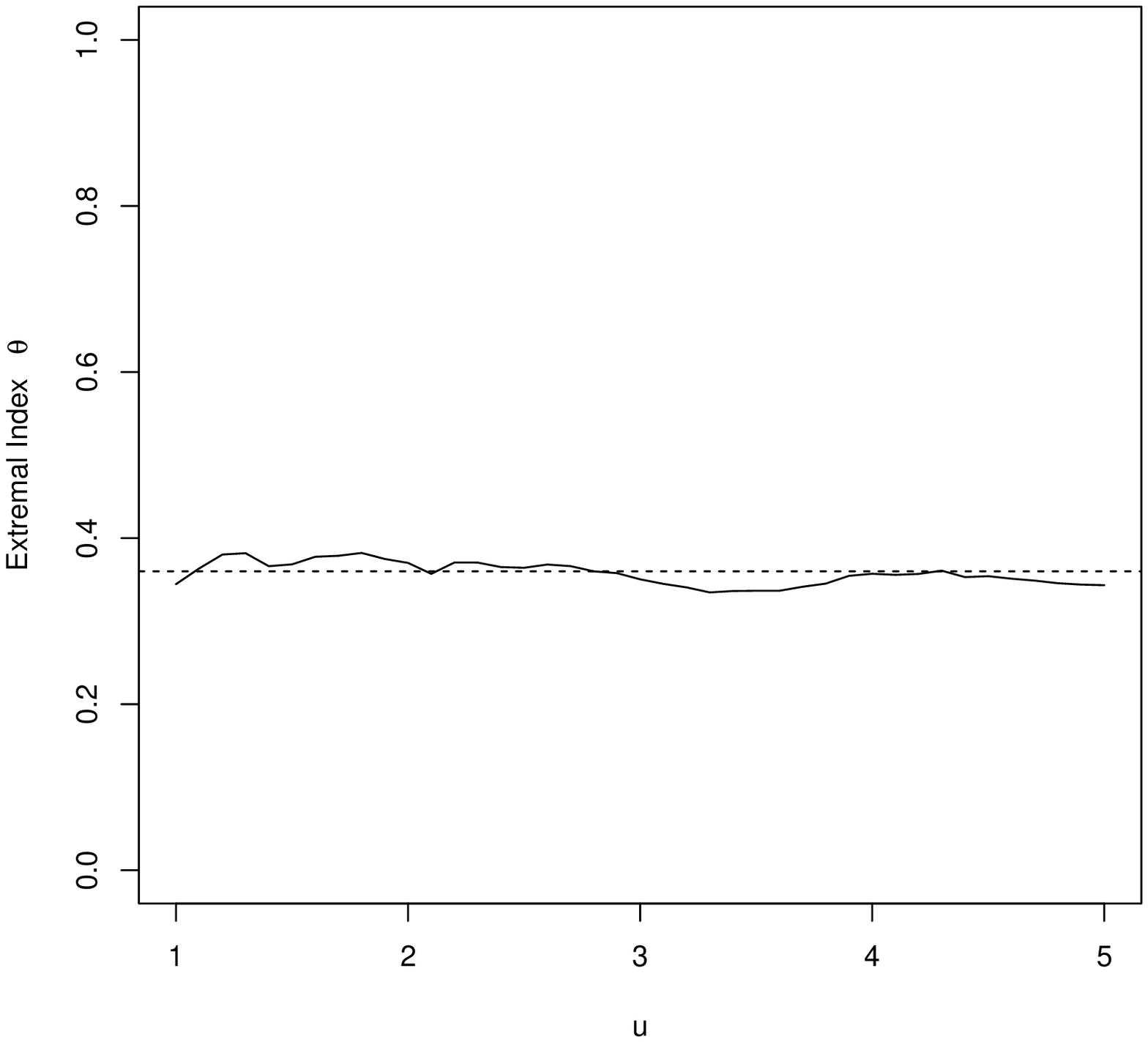,height=6cm,width=6cm,angle=0}}
\bfi{\em  Blocks estimator $\wh
  \theta$ of the extremal index $\theta$ for a sample of size  $20\,000$ from
the  AR(1) process $X_t=0.8\,X_{t-1}+Z_t$. The iid noise 
$(Z_t)$ has a common student \ds\ with $\alpha=2$ degrees of freedom.
The extremal index $\theta=0.37$ is known (indicated by dashed line); see 
\cite{embrechts:kluppelberg:mikosch:1997}, Section 8.1.
The blocks
  estimator as a \fct\ of the block size $s$ and 
$u\%$ of the upper order statistics
(top), for fixed $u=2\% $ and running $s$ (bottom left) and for fixed
$s=24$ and running $u$ (bottom right). 
 }\label{fig:index2}
\efi
\end{figure}

In Figures \ref{fig:index1}
and \ref{fig:index2} we illustrate the estimation of $\theta$ for real
and simulated data. We choose the simple {\em blocks estimator} $\wh
  \theta=K_n/N_n$ of the extremal index $\theta$,
where $N_n$ is the number of exceedances of the threshold $u=u_n$ in the
sample $X_1,\ldots,X_n$ and $K_n$ is the number of blocks of size $s=s_n$,
$X_{(i-1)s+1},\ldots,X_{is}$, $i=1,\ldots,[n/s]$, with at least one
exceedance of $u$.

Aspects of bias, variance
and optimal choice of blocks for the estimation of $\theta$ were 
discussed in Smith and Weissman \cite{smith:weissman:1994}. In a
series of papers, Hsing \cite{hsing:1988,hsing:1991a,hsing:1991b,hsing:1993}
studied the extremes of stationary \seq s,
including  the \asy\ behavior of their extremal index estimators.
The recent papers Robert \cite{robert:2009a,robert:2009b},
Robert et al. \cite{robert:segers:ferro:2009}, in particular 
\cite{robert:2009b}, give historical accounts of estimation of
$\theta$ and some new technology for the estimation of $\theta$ 
and the cluster size \ds\ $P(\xi=k)$, $k\ge 1$. The paper of Robert
\cite{robert:2009b} is devoted to inference on the cluster size \ds .
The literature on this topic is sparse; Robert \cite{robert:2009b} mentions Hsing 
\cite{hsing:1991a} as a historical reference.
\subsection{The extremogram: an \asy\ correlogram for extreme events}\label{subsec:extremogram}
Davis and Mikosch \cite{davis:mikosch:2009a,davis:mikosch:2012}
introduced another tool for measuring the extremal dependence in a 
strictly stationary $\bbr^d$-valued \ts\ $(X_t)$: 
the {\em extremogram} defined as a
limiting \seq\ given by
\beam\label{eq:extrem}
\gamma_{AB}(h)=\lim_{\nto}n\,\cov(I_{\{a_n^{-1} X_0\in A\}} ,I_{\{a_n^{-1} X_h\in B\}})\,,\quad
  h\ge 0\,.
\eeam
Here $(a_n)$ is a suitably chosen normalization \seq\
and $A,B$ are two fixed sets bounded away from zero. 
The events $\{X_0\in a_n A\}$ and 
$\{X_h\in a_n B\}$ are considered as extreme ones and $\gamma_{AB}(h)$
measures the influence of the time zero extremal event
$\{X_0\in a_n\,A\}$ on the extremal event $\{X_h\in
a_n\,B\}$, $h$ lags apart.
The choice of $(a_n)$ depends on the situation at hand.   
To avoid ambiguity, we later  assume that
$(a_n)$ satisfies the relation $n\, P(|X|>a_n)\sim 1$.
With this choice of $(a_n)$, $\gamma_{AB}(h)=
\lim_{\nto}n\,P(a_n^{-1}X_0\in A,a_n^{-1} X_h\in B)$. 
Motivating examples of extremograms are the limiting conditional 
\pro ies
$\lim_{\nto} P(a_n^{-1}X_h\in B\mid a_n^{-1} X_0\in A)$
in  Davis and Mikosch \cite{davis:mikosch:2009a,davis:mikosch:2012}. 
\par
A motivating example for $d=1$
with $A=B=(1,\infty)$ is the so-called {\em (upper)
tail dependence coefficient} of the vector $(X_0,X_h)$ given as the limit
\beam\label{eq:rho}
\rho(h)=\lim_{\xto} P(X_h>x\mid X_0>x)\,.
\eeam
(Here we assume that $X$ has infinite right endpoint.) These pairwise tail
dependence coefficients have attracted some attention in the
literature on quantitative risk management; see for example 
McNeil et al. \cite{mcneil:frey:embrechts:2005}.
Notice that
$\rho(h)$ coincides with $\gamma_{AA}(h)$ if we  choose $(a_n)$ \st\
$n P(X_0>a_n)\sim 1$ as $\nto$. Indeed, 
\beao
n \, \cov(I_{\{X_0>a_n\}} ,I_{\{X_h>a_n
  \}})
&\sim & \dfrac{P(X_h>a_n,X_0>a_n)-(P(X_0>a_n))^2}{P(X_0>a_n)}\\
&\sim & P(X_h>a_n\mid X_0>a_n)\,.
\eeao
\par
A similar calculation for any dimension $d$ and suitable sets $A,B$
shows that the limiting \seq
\beam\label{eq:kk}
\Big(\barr{cc}\gamma_{AA}(h)&\gamma_{AB}(h)\\
\gamma_{BA}(h)&\gamma_{BB}(h)
\earr\Big)\,, \quad h\ge 0\,,
\eeam 
inherits the properties of a matrix covariance \fct . Notice that
the entries of these matrices cannot be negative. The 
interpretation of \eqref{eq:kk} as covariance \fct\ allows one to 
use the classical notions of \tsa\ {\em in an \asy\ sense.} 
For example, notions such as long or short range dependence of
extremal events can be made precise by specifying the rate of decay of
\eqref{eq:kk} as $h\to \infty$. Davis and Mikosch
\cite{davis:mikosch:2009c},
Mikosch and Zhao \cite{mikosch:zhao:2012}
introduced an analog of the spectral density as a Fourier transform of
the \seq\ \eqref{eq:kk}. They showed that the \per\ of the 
\seq\ of indicators $I_{\{a_n^{-1} X_t\in A\}}$ of the
extremal events $\{a_n^{-1} X_t\in A\}$, $t\in\bbz$, has properties
similar to the classical \per\ of a stationary \seq .  In particular,
weighted averages of the \per\ are consistent estimators of the 
spectral density.
\par
In the literature,  the pairwise tail dependence 
coefficients \eqref{eq:rho} are 
mostly considered for concrete examples of \ds s, such as elliptical
ones, including the multivariate $t$- and Gaussian \ds s; see e.g.
McNeil et al. \cite{mcneil:frey:embrechts:2005}. In these
cases, one can verify that the limits in \eqref{eq:rho} exist. 
For a Gaussian stationary \seq , $\rho(h)=0$, $h\ge 1$, unless $X_t=X$
a.s.  for all $t\in\bbz$. The case $\rho(h)=0$ for some $h\ge 1$ is
(again) referred to as {\em \asy\ extremal dependence} in the vector 
$(X_0,X_h)$ although no extremal index is in view.
\par
In
general, it is not obvious whether the limits $\rho(h)$ and, more
generally, $\gamma_{AB}(h)$ for $h\ge 0$ exist. In this paper, we will
use the notion of a {\em \regvary\ stationary \seq }. It is a
sufficient condition for the existence of the limits $\gamma_{AB}(h)$.
Roughly speaking, a \regvary\ \seq\ of \rv s $(X_t)$ has power law
tails for every lagged vector $(X_1,\ldots,X_h)$, $h\ge 1$. 
In what follows,  we make precise what \regvar\ means.

\subsubsection{Regularly varying random vectors}\label{subsec:rv}
The notion of \regvar\ is basic in \evt\ and limit theory for
partial sums of iid \rv s. In multivariate \evt , \regvar\ with index
$\alpha>0$ of the 
$d$-dimensional iid random vectors $X_t$, $t\in\bbr$,
with values in  $(0,\infty)^d$ is necessary and sufficient for the fact 
 that the normalized 
\seq\ of component-wise maxima 
$(a_n^{-1} \max_{t\le n} X_t^{(i)})_{i=1,\ldots,d}$, $t=1,2,\ldots,$
converges in \ds\ to a $d$-dimensional \evd\ $H$ on  $(0,\infty)^d$
whose marginal \ds s are Fr\'echet
$\Phi_\alpha$-distributed; see Resnick \cite{resnick:1987} for a
general theory of multivariate extremes for iid \seq s.
 Similarly, for a general $\bbr^d$-valued iid  
\seq\ $(X_t)$, the \seq\ of suitably normalized and centered partial sums 
$a_n^{-1}(X_1+\cdots+ X_n-b_n)$ converges in \ds\ to an infinite
variance $\alpha$-stable limit \fif\ the \ds\ of $X$ is \regvary\ with
index $\alpha$. The index $\alpha$ is then necessarily in the range
$\alpha\in (0,2)$. We refer to Rva\v ceva \cite{rvaceva:1962} and
Resnick \cite{resnick:2007}
for proofs of this fact.
\par
Various definitions of a $d$-dimensional \regvary\ vector $X$ 
exist; we refer to Resnick \cite{resnick:1986,resnick:1987,resnick:2007}.
 We start with a definition in terms of spherical coordinates.
We say that $X$ is  \regvary\ with index $\alpha>0$ and spectral
\ms\ $P(\Theta \in \cdot)$ on the Borel $\sigma$-field of 
the unit sphere 
$\bbs^{d-1}=\{x\in\bbr^d: |x|=1\}$ if \footnote{The choice of the norm
  $|\cdot|$ is relevant for defining the corresponding unit sphere and
  the spectral \ms\ on it, but the notion of \regvar\ of a vector does
  not depend on a particular choice of norm. In this paper, $|\cdot|$
  will stand for the Euclidean norm.}
the following weak limits exist
for every fixed $t>0$:
\beam\label{eq:spec}
\dfrac{P(|X|>tx\,, X/|X|\in \cdot)}{P(|X|>x)}\stw t^{-\alpha}
\,P(\Theta \in \cdot)\,,\quad \xto\,.
\eeam    
Relation \eqref{eq:spec} can be written in an equivalent form as
a pair of conditions: 
\begin{enumerate}
\item
The norm $|X|$ is \regvary\ in the classical sense,
i.e. $P(|X|>tx)/P(|X|>x)\to t^{-\alpha}$, $t>0$, or, equivalently,
$P(|X|>x)= x^{-\alpha}L(x)$, $x>0$, for a \slvary\ \fct\ $L$; cf. 
Bingham et al. \cite{bingham:goldie:teugels:1987} for an encyclopedia
on \regvary\ \fct s.
\item
The angular component  $X/|X|$ is independent of $|X|$ for large
values of $|X|$ in the sense that
\beam\label{eq:spec1}
P(X/|X|\in \cdot \mid |X|>x)\stw P(\Theta\in \cdot)\,,\quad \xto\,.
\eeam
\end{enumerate}
In any of these limit relations, it is possible to replace the
converging parameter $x$ by a \seq\ $(a_n)$ \st\ $P(|X|>a_n)\sim
n^{-1}$.
Then \eqref{eq:spec} and \eqref{eq:spec1}, respectively, read as 
\beao
n\,P(|X|>t a_n, X/|X|\in \cdot)\stw t^{-\alpha} P( \Theta\in \cdot)\quad
\mbox{and}\quad  P(X/|X|\in \cdot \mid |X|>a_n)\stw P(\Theta\in \cdot)\,.
\eeao
The \con\ relation \eqref{eq:spec} can be understood as \con\ on the 
particular Borel sets $\{\bfx\in \bbr^d: |\bfx|>t, \bfx/|\bfx|\in S\}$
for Borel sets $S\subset \bbs^{d-1}$ with a smooth boundary. This \con\ can be extended to 
the Borel $\sigma$-field on $\ov \bbr_0^d=\ov \bbr^d \setminus
\{\bf0\}$, $\ov \bbr = \bbr \cup \{\infty,-\infty\}$:
\beam\label{eq:spec2}
\mu_x(\cdot)=\dfrac{P(x^{-1 } X\in \cdot) }{P(|X|>x)} \stv \mu(\cdot)\,,\quad \xto\,.
\eeam
Here $\stv$ refers to vague \con\ of \ms s on 
the Borel $\sigma$-field on $\ov \bbr_0^d$, i.e. $\int_{\ov \bbr_0^d} f\,
d\mu_x\to \int_{\ov \bbr_0^d} f\,
d\mu$ as $\xto$ for any continuous and compactly supported $f$ on $\ov
\bbr_0^d$; see Kallenberg \cite{kallenberg:1983}, Resnick \cite{resnick:1987}.
This means in particular, that the support of $f$ is bounded away from
zero. In view of \eqref{eq:spec}, $\mu (\{\bfx\in \ov \bbr^d_0: |\bfx|>t,
\bfx/|\bfx|\in S\})= t^{-\alpha} P(\Theta\in S)$, and therefore $\mu$
is a Radon \ms\ (i.e. finite on sets bounded away from zero) satisfying
$\mu (t A)=t^{-\alpha} \mu(A)$, $t>0$. In particular, $\mu$ does not charge
points containing infinite components. 
Again, the parameter $x$
in \eqref{eq:spec2} can be replaced by a \seq\ $(a_n)$ satisfying 
$P(|X|>a_n)\sim n^{-1}$ and then we get
\beao
n\,P(a_n^{-1}X\in \cdot )\stv \mu(\cdot) \,,\quad \nto\,.
\eeao
For an iid \seq\ $(X_t)$ with generic element $X$, the latter condition
is equivalent to the \con\ of the \pp es
\beao
N_n= \sum_{t=1}^n \vep_{a_n^{-1} X_t}\std N\,,
\eeao
where $N$ is a Poisson random \ms\ with mean \ms\ $\mu$ and state space
$\ov \bbr_0^d$; see Resnick \cite{resnick:1986,resnick:1987}. Since \pp\ \con\ is
basic to \evt , the  notion of multivariate \regvar\ is very natural
in the context of \evt\ for multivariate observations with
heavy-tailed components; see also the recent monograph by Resnick 
\cite{resnick:2007} who stresses the importance of the notion of
\regvar\ as relevant for many applications in finance, \ins\ 
and telecommunications.\footnote{The notion of \regvar\ is essentially
dimensionless; see for example relation \eqref{eq:spec} which immediately 
extends to normed spaces and, more generally, to metric spaces. An
account of the corresponding theory can be found in Hult and Lindskog 
\cite{hult:lindskog:2006}. Applications of  \regvar\ 
in \fct\ spaces to \evt\ can be found in 
de Haan and Tao \cite{haan:tao:2003},  
Davis and Mikosch \cite{davis:mikosch:2008}, 
Meinguet and Segers \cite{meinguet:segers:2010},
to \ld s in Hult et al. 
\cite{hult:lindskog:mikosch:samorodnitsky:2005}, Mikosch and
Wintenberger \cite{mikosch:wintenberger:2012a,mikosch:wintenberger:2012b}, 
and to random sets in 
Mikosch
et al.
\cite{mikosch:pawlas:samorodnitsky:2011a,mikosch:pawlas:samorodnitsky:2011b}.}
\subsubsection{Regularly varying stationary \seq
  s}\label{subsec:regvar}
A strictly stationary \seq\ $(X_t)$ is \regvary\ with index $\alpha$
if its \fidi s are \regvary\ with index $\alpha$, i.e. for every $h\ge
1$, there exist 
non-null Radon \ms s $\mu_h$ on the Borel $\sigma$-field of $\ov
\bbr_0^h$ and a \seq\ $(a_n)$  \st\ $a_n\to\infty$ and 
\beam\label{eq:muh}
n\,P(a_n^{-1}(X_1,\ldots,X_h)\in \cdot) \stv \mu_h(\cdot)\,,\quad \nto\,.
\eeam
Here and in what follows, we will choose the normalizing \seq\ $(a_n)$
\st\  $P(|X|>a_n)\sim n^{-1}$, where we use the
notation $(a_n)$  in a way different from
Section~\ref{subsec:rv}. Indeed, in the latter section we defined
the normalization $(a_n^{(h)})$ \st\ 
$P(|(X_1,\ldots,X_h)|>a_n^{(h)})\sim n^{-1}$, but then the
normalization would depend on the dimension $h$. This is not
desirable.
However, notice that 
\beao
1&=& \lim_{\nto}\dfrac{P(|(X_1,\ldots,X_h)|>a_n^{(h)})}{P(|X|>a_n)}\,.
\eeao
Therefore, by the properties of \regvary\ \fct s, there exist
positive constants $c_h^{1/\alpha}=\lim_{\nto} a_n/a_n^{(h)}$, $h\ge 1$. Hence
for sets $A$ bounded away from zero \st\  $\mu_h(\partial A)=0$ we have
\beao
n\,P((a_n^{(h)})^{-1}(X_1,\ldots,X_h)\in A)&=&n \, 
P(a_n^{-1}(a_n/a_n^{(h)})(X_1,\ldots,X_h)\in A)\\
&\sim& c_h\, [n\,P(a_n^{-1}(X_1,\ldots,X_h)\in A)] \\
&\to & c_h\, \mu_h(A)\,,
\eeao
i.e. the limit \ms s of \regvar\ under the different normalizations only
differ by some positive constants.
\par
The condition of \regvar\ on the \seq\ $(X_t)$ seems to be a severe
restriction since the tails of
the marginals are power laws. 
However, following Resnick \cite{resnick:1987}, Proposition~5.10, any
multivariate \ds\ (with continuous marginals) 
in the maximum domain of attraction (MDA) 
of a $d$-dimensional \evd\ can be transformed
to a \ds\ $G$ with common Fr\'echet or Pareto marginals. 
Then $G$ is in the MDA 
of an \evd\ with Fr\'echet marginals or,
equivalently, $G$ is \regvary .
\par
For example, transforming the marginals of a Gaussian stationary \seq\
to unit Fr\'echet, the resulting \seq\ is \regvary\ with index
$\alpha=1$.
We mentioned before that the tail dependence coefficient $\rho(h)=0$,
$h\ge 1$, for
any non-trivial Gaussian stationary \seq . The quantities
$\rho(h)$ remain invariant under monotone 
increasing transformations of the 
marginals. Hence, the transformed Gaussian \ds\ with
unit Fr\'echet marginals exhibits \asy\ independence in the sense that
the limit \ms s $\mu_h$ are concentrated on the axes.
\subsubsection*{The tail process} 
An insightful characterization of an $\bbr^d$-valued \regvary\
stationary \seq\  $(X_t)$ was given in
Theorem 2.1 of Basrak and Segers \cite{basrak:segers:2009}: 
there exists a \seq\ of $\bbr^d$-valued random vectors
$(Y_t)_{t\in\bbz}$ \st\ $P(|Y_0| > y) = y^{-\alpha}$ for $y > 1$ and
for any $h\ge 0$,
\beao
P (x^{-1}(X_{-h},\ldots,X_h)\in\cdot \mid |X_0| > x)\stw
P((Y_{-h},\ldots,Y_h)\in \cdot)\,,\quad \xto\,.
\eeao
The process
$(Y_t)$ is the {\em tail process} of $(X_t)$. 
Writing $\Theta_t = Y_t/|Y_0|$ for $t\in \bbz$,
one also has for $h\ge 0$,
\beao
P( |X_0|^{-1}(X_{-h},\ldots,X_h)\in\cdot \mid |X_0| > x)
\stw P ((\Theta_{-h},\ldots,\Theta_h)\in \cdot)\,,\quad \xto\,.
\eeao
The process $(\Theta_t)$ is independent of $|Y_0|$ and  called
the {\em spectral tail process} of $(X_t)$. Notice that
$P(\Theta_0\in \cdot)$ is the spectral \ms\ of $X$.
\par
Basrak and Segers \cite{basrak:segers:2009} also gave an expression
for the extremal index in terms of the spectral tail process:
\beam\label{eq:basrak}
\theta= E\big[\sup_{t\ge 0} |\Theta_t|^\alpha-\sup_{t\ge 1} 
|\Theta_t|^\alpha\big]\,.
\eeam

\subsubsection{The extremogram revisited}
Now consider an $\bbr^d$-valued \regvary\ stationary $(X_t)$. 
Then the extremogram $\gamma_{AB}(h)$, $h\ge 0$, is well defined. Indeed,
for every $h\ge 0$,  the vector $(X_1,\ldots,X_{h+1})$ is \regvary\ with
limit \ms\ $\mu_{h+1}$. Then, with normalization $(a_n)$
\st\ $P(|X|>a_n)\sim n^{-1}$,
\beao
n\,P(a_n^{-1} X_{h}\in B, a_n^{-1}X_0 \in A)\to \mu_{h+1}(A\times
\bbr^{d (h-1)}\times B)=\gamma_{AB}(h)\,, \quad h\ge 0\,,
\eeao  
provided 
$A\times \bbr^{d (h-1)}\times B$ is a continuity set \wrt\ 
the \ms\ $\mu_{h+1}$. Similarly, for $d=1$ and $A=B=(1, \infty)$,
\beao
\rho(h) = \dfrac{\mu_{h+1}(A\times \bbr^{h-1}\times  A)}{\mu_{h+1}(A\times
  \bbr^h)}\,,\quad h\ge 0\,.
\eeao

These limits can also be expressed in terms of the tail process.
In the former case, assuming that $A$ is bounded away from zero, there
exists $\delta>0$ \st\ $A\subset \{x\in\bbr^d: |x|>\delta \}$. Hence 
\beao
\lefteqn{\dfrac{P(a_n^{-1} X_{h}\in B, a_n^{-1}X_0 \in A)}{P(|X|>a_n)}}\\
&=&\dfrac{P(a_n^{-1} X_{h}\in B, a_n^{-1}X_0 \in A, |X_0|>\delta
  a_n)}{P(|X|>a_n)}\\
&=&\dfrac{P((\delta a_n)^{-1} X_{h}\in \delta^{-1}B, (\delta a_n)^{-1}X_0 \in \delta^{-1}A, |X_0|>\delta
  a_n)}{P(|X|>\delta a_n)}\dfrac{P(|X|>\delta a_n)}{P(|X|>a_n)}\\
&\to &P( (Y_0,Y_h)\in \delta^{-1} (A\times B)) \,\delta^{-\alpha}\\&=&
P( (Y_0,Y_h)\in  A\times B)=\gamma_{AB}(h)\,. 
\eeao
Similarly, for $d=1$ and $A=B=(1,\infty)$, assuming that 
$\lim_{\xto} P(X>x)/P(|X|>x)=E(\Theta_0)_+^\alpha=P(\Theta_0=1)>0$,
\beao
\dfrac{P(X_{h}>a_n, X_0 >a_n)}{P(X>a_n)}
&=&\dfrac{P(X_{h}>a_n, X_0 >a_n, |X_0|>a_n)}{P(|X|>a_n)}
\dfrac{P(|X|>a_n)}{P(X>a_n)}\\
& \to & P(Y_h>1\mid Y_0>1)\\&=&\dfrac{P(|Y_0|
  \min(\Theta_0,\Theta_h)>1)}{P(|Y_0| \Theta_0>1)}\\&=&\dfrac{E(\min(\Theta_0, \Theta_h))_+^\alpha}{E(\Theta_0)_+^\alpha}=\rho(h)\,.
\eeao
\subsubsection{Examples of \regvary\ \seq s and their extremograms}\label{subsub:examples} 
In this section, we will introduce some important classes of
real-valued strictly stationary \regvary\
stationary  \seq s with index $\alpha>0$. We will also give the values
of the extremogram $\rho(h)$, $h\ge 1$, in \eqref{eq:rho}.
For the calculation of $\rho$ in these examples, we refer to 
\cite{davis:mikosch:2009c,davis:mikosch:2012,mikosch:zhao:2012}.
\subsubsection*{IID \seq } 
An iid \seq\ $(Z_t)$ is \regvary\ with index $\alpha$ 
\fif\ $Z$ is \regvary\ with the same index;  the limit \ms s $\mu_h$ are concentrated on
the axes
and $\rho(h)=0$, $h\ge 1$.
\subsubsection*{Linear process}
Historically, the class of linear processes with
\regvary\ iid real-valued noise $(Z_t)$ has attracted attention in \evt\ and in
\tsa . A (causal) linear process
\beam\label{eq:lin}
X_t = \sum_{t=0}^\infty \psi_j Z_{t-j}\,,\quad t\in\bbz\,,
\eeam
inherits \regvar\ under conditions on the deterministic \seq\
$(\psi_i)$ which are close to those dictated by the 3-series theorem,
ensuring the a.s. \con\ of the series in \eqref{eq:lin}.
This fact was proved in Mikosch and Samorodnitsky 
\cite{mikosch:samorodnitsky:2000} for the \ds\ of $X$. The \regvar\ of
the \fidi s of $(X_t)$ follows since \regvar\ is preserved under
affine transformations of \regvary\ vectors. The class \eqref{eq:lin}
includes causal ARMA processes which are relevant for
applications. We refer to Chapter~7 of Embrechts et al. 
\cite{embrechts:kluppelberg:mikosch:1997} for various applications of 
\regvary\ linear processes. 
\par
Under the tail balance condition $P(Z>x)\sim p\,P(|Z|>x)$, 
$P(Z\le -x)\sim q\,P(|Z|>x)$, as $\xto$, for some $p,q\ge 0$ with
$p+q=1$,
\beao
\rho(h) = \dfrac{\sum_{i=0}^\infty \Big[ p\, (\min (\psi_{i}^+ ,  
  \psi_{i+h}^+))^\alpha + q\, (\min (\psi_{i}^- ,  
  \psi_{i+h}^-))^\alpha \Big]}{ \sum_{i=0}^{\infty} \Big[ p\,  
  (\psi_i^{+})^\alpha+ q (\psi_i^{-})^\alpha\Big]}\,,\quad h\ge 1\,.  
\eeao

\subsubsection*{Stochastic recurrence equations} Next to linear processes,
solutions to the stochastic recurrence equation
\beam\label{eq:sre}
X_t=A_t\,X_{t-1}+B_t\,,\quad t\in\bbz\,,
\eeam
have attracted some attention. Here $(A_t,B_t)$, $t\in\bbz$, is an iid
$\bbr^2$-valued \seq . An a.s unique causal solution to \eqref{eq:sre}
exists under the moment conditions $E\log A^+<0$ and $E\log
|B|<\infty$. It follows from work by Kesten \cite{kesten:1973} and
Goldie \cite{kesten:1973} that $X$ is \regvary\ in the precise sense
that
\beao
P(X>x)\sim c_+ x^{-\alpha}\quad \mbox{and}\quad 
P(X \le -x)\sim c_- x^{-\alpha}\,,\quad \xto \,, 
\eeao
for constants $c_+,c_-\ge 0$ \st\ $c_++c_->0$ provided the equation
\beam\label{eq:kesten}
E|A|^\kappa=1
\eeam 
has a positive solution $\alpha$ (which is unique due
to convexity), $EB^\alpha<\infty$ and further regularity conditions
on the \ds\ of $A$ are satisfied. Iteration of \eqref{eq:sre} shows
that the \fidi s of $(X_t)$ are \regvary\ with index $\alpha$.
This fact is rather surprising since the \ds s of $A$ and $B$ do not
need to be heavy-tailed, in contrast to linear processes , where the
noise $(Z_t)$ itself has to be heavy-tailed to ensure  \regvar\ of
$(X_t)$. We mention that the case of multivariate $B$ and
matrix-valued $A$ has also been studied, starting with Kesten
\cite{kesten:1973}; see the recent paper 
Buraczewski et al. \cite{buraczewski:damek:guivarch:hulanicki:2009}.
\par
Assuming $A>0$ a.s., similar calculations as in the proof of Lemma 2.1 in 
\cite{davis:mikosch:2009c} yield
\beam\label{eq:kka}
\rho(h)= E[\min(1,A_1\cdots A_h)^\alpha]\,,\quad h\ge 1\,.
\eeam
\subsubsection*{Models for returns} Log-returns $X_t=\log P_t-\log
P_{t-1}$, $t\in\bbz$,  of a speculative price series $(P_t)$ 
are often modeled of the form $X_t=\sigma_t Z_t$, where $(\sigma_t)$
is a strictly stationary \seq\ of non-negative {\em volatilities} and
$(Z_t)$ is an iid multiplicative noise \seq . The feedback between
$(\sigma_t)$ and $(Z_t)$ can be modeled in a rather flexible way.
\subsubsection*{Stochastic volatility models}
The most simple approach is to assume that $(\sigma_t)$ and $(Z_t)$ be
independent.  The resulting time series model is frequently referred to as 
{\em stochastic volatility model}. Its probabilistic properties are
rather simple; see Davis and Mikosch \cite{davis:mikosch:2009d}.
In particular, \regvar\ of $(X_t)$ results if 
$E\sigma^{\alpha+\delta}<\infty$ for some $\delta>0$ and $(Z_t)$ is
iid and \regvary\ with index $\alpha$. The corresponding limit \ms s
$\mu_h$ in \eqref{eq:muh} are then concentrated on the axes; see Davis
and 
Mikosch \cite{davis:mikosch:2001,davis:mikosch:2009b},\footnote{The fact that $\mu_h$, $h\ge 1$, is concentrated on the axes is also
referred to as {\em \asy\ extremal independence.} Recall that, in an extreme value
context, various other situations are also referred to as \asy\
extremal independence, among them the cases of unit extremal index and
zero tail dependence coefficient. Asymptotic independence in
the sense of the limiting \ms s $\mu_h$ is much more complex than the other 
notions which are just numerical characteristics. The fact that
$\mu_h$ is concentrated on the axes means that it is very unlikely
that any two values $X_t$ and $X_s$, $s\ne t$, are big at the same
time, just as for independent \rv s. On the other hand, this kind of \asy\ 
independence heavily relies on the notion of multivariate \regvar .}
and then also $\rho(h)=0$, $h\ge 1$, as in the iid case.
The situation changes if $E|Z|^{\alpha+\delta}<\infty$ for some
$\delta>0$ and $(\sigma_t)$ is \regvary\ with index $\alpha$. Then 
$(X_t)$ is \regvary\ with index $\alpha$ and extremal clustering for
this \seq\ is possible; see Mikosch and Rezapur
\cite{mikosch:rezapur:2012}. 

\subsubsection*{GARCH model}
Among the models for returns $X_t=\sigma_t Z_t$, $t\in\bbz$, 
the GARCH family gained most popularity.
The simplest model of its kind (ARCH) was introduced by Engle 
\cite{engle:1982} and the more sophisticated GARCH model by Bollerslev
\cite{bollerslev:1986}.
For simplicity, we consider
the \garch\ case given by $\sigma_t^2=\alpha_0 + \sigma_{t-1}^2
(\alpha_1Z_{t-1}^2+\beta_1)$, $t\in\bbz$, where $\alpha_0,\alpha_1,
\beta_1$ are positive constants with certain restrictions on the
values of $\alpha_1+ \beta_1,\beta_1<1,$ to ensure strict stationarity. Typical
choices are  standard normal or unit variance $t$-distributed
$Z$. Notice that we can write $\sigma_t^2= A_t\sigma_{t-1}^2+B_t$,
where $B_t=\alpha_0$ and $A_t=\alpha_1Z_{t-1}^2+\beta_1$, $t\in\bbz$.
Therefore \regvar\ of $(\sigma_t)$ follows from the corresponding
theory for \sre s, and $(X_t)$ inherits the same property; see 
Davis and Mikosch \cite{davis:mikosch:1998} for the ARCH(1) case,
Mikosch and \sta\ \cite{mikosch:starica:2000} for the \garch ,
Basrak et al. \cite{basrak:davis:mikosch:2002} for \garchpq\ and the
review paper Davis and Mikosch \cite{davis:mikosch:2009a}. 
Real-life log-returns are typically heavy-tailed. The GARCH model
captures this property and this was one of the reasons that it became 
a benchmark model in 
financial \tsa\ from which numerous other models where derived. An
expression
of $\rho$ for $(X_t^2)$ is given by \eqref{eq:kka} with $A_t=\alpha_1 Z_t^2+\beta_1$. 

\subsubsection*{Infinite variance stable \seq } Stable processes with infinite
variance have become popular due to their attractive theoretical and
modeling properties; see Samorodnitsky and Taqqu 
\cite{samorodnitsky:taqqu:1994}. The \fidi s of an $\alpha$-stable process are
jointly $\alpha$-stable, hence they are \regvary\ with 
index $\alpha\in (0,2)$.  The class of infinite variance
stationary stable processes has been intensively studied; see 
Rosi\'nski \cite{rosinski:1995}. An expression of $\rho$ is given in 
\cite{davis:mikosch:2009c}, Section 2.4.

\subsubsection*{Max-stable processes} This class of processes has
recently attracted some attention since it is a flexible class for modeling 
heavy tails and spatio-temporal dependence. Since the \fidi s of
max-stable processes are explicitly given it is often simple to verify
properties (such as \regvar ) and to calculate
certain quantities (e.g. mixing coefficients, extremal index). 
We will use this class of \regvary\ processes to illustrate the general theory.
\par
Following de Haan \cite{haan:1984}, a real-valued process
$(\xi_t)_{t\in T}$, $T\subset \bbr$, is {\em
  $\alpha$-max-stable} for some $\alpha>0$ if its \fidi s satisfy the relation
\beao
P(\xi_{t_1}\le x_1,\ldots, \xi_{t_d}\le x_n)&=&\exp
\Big\{ - \int _{\bbs^{d-1}\cap \bbr_+^d} 
\max_{i\le d} \big(\dfrac{s_i}{x_i}\big)^\alpha\,  \Gamma_ {\bft_d}(d\bfs )\Big\}\,,\\&&\quad t_i\in
T\,,i=1,2,\ldots,d\,, \quad x_i>0\,,\quad d\ge 1\,,
\eeao  
where $\Gamma_ {\bft_d}$ are finite \ms s on the unit sphere.
This means in particular that the marginal \ds s of the process $\xi$
have a Fr\'echet  \ds\ with parameter $\alpha$ given by\footnote{The choice of
Frech\'et $\Phi_\alpha$ marginals is for convenience only; then 
results on \regvar\ are applicable.  Since 
Gumbel or Weibull distributed \rv s can be obtained by suitable increasing 
transformations of
a Fr\'echet \rv\ any result for max-stable processes with Fr\'echet
marginals  can be formulated in terms of the transformed processes
with Gumbel or Weibull marginals; see for example Kabluchko et. al
\cite{kabluchko:schlather:haan:2009} who formulated their results
in terms of Gumbel \ds s. Since the choice of the parameter $\alpha$
is also arbitrary in this context, most results in the literature 
are formulated for processes with unit Fr\'echet $\Phi_1$ marginals.}
\beam\label{eq:frechet}
\Phi_\alpha(x)= \ex^{-x^{-\alpha}}\,,\quad x>0\,.
\eeam
 De Haan
\cite{haan:1984}
also introduced the notion of {\em $\alpha$-max-stable
integral}. Given a $\sigma$-finite \ms\ space $(E,{\mathcal E},\nu)$,
consider a Poisson random \ms\ $\sum_{i=1}^\infty
\vep_{(\Gamma_i,Y_i)}$
with $0<\Gamma_1<\Gamma_2<\cdots$
on the state space $\bbr_+\times E$
with mean \ms\ $\Leb \times \nu$.  For $f\ge 0$
with $f\in L^\alpha(E,{\mathcal E},\nu)$ the max-stable integral is
defined as 
\beao
\int_E^{\vee} f\, d M_\nu^\alpha= \sup_{i\ge 1} 
\Gamma_i^{-1/\alpha}\,f(Y_i)\,\,. 
\eeao 
Using the  order statistics property of the homogeneous Poisson
process with points $(\Gamma_i)$, one obtains 
\beam\label{eq:maxs}
P \big(\int_E^{\vee} f\, d M_\nu^\alpha\le x\big)=
\exp\Big\{-x^{-\alpha}\int_E f^\alpha(y)\, \nu(dy)\Big\}\,,\quad x> 0\,.
\eeam
Moreover, for any non-negative  $f_i\in L^\alpha(E,{\mathcal E},\nu)$,
$i=1,\ldots,d$, by \eqref{eq:maxs},
\beam
P\Big(\int_E^{\vee} f_i \,d M_\nu^\alpha\le x_i\,,i=1,\ldots,d\Big)
&=&P\Big(\int_E \max_{i=1,\ldots,d}\dfrac{f_i}{x_i}\, d M_\nu^\alpha\le
1\Big)\nonumber\\
&=&\exp\Big\{-\int_E \max_{i=1,\ldots,d}\big(\dfrac{f_i(y)}{x_i}\big)^\alpha
\nu(dy)\Big\}\,,\quad x_i>0\,, i=1,\ldots,d\,.\label{eq:fidi}
\eeam
The notions of max-stable process and integral bear some
resemblance with the corresponding $\alpha$-stable ones; see Stoev
and Taqqu \cite{stoev:taqqu:2005}, Kabluchko \cite{kabluchko:2009}.
\par
We will focus on {\em stationary ergodic} max-stable processes with
integral \rep\ 
\beam\label{eq:repr}
X_t=\int_E^{\vee} f_t\, d M_\nu^\alpha\,,\quad t \in \bbr\,,\quad f_t\ge 0\,,
\quad f_t\in L^\alpha(E,{\mathcal E},\nu)\,.
\eeam
As in the case of $\alpha$-stable stationary
ergodic processes (Rosi\'nski \cite{rosinski:1995}), 
the choice of $(f_t)$ is 
rather sophisticated; see Stoev \cite{stoev:2008}, Kabluchko
\cite{kabluchko:2009} for details. De Haan \cite{haan:1984}
showed that any max-stable process with countable index set $T\subset \bbr$ 
and stochastically continuous sample paths has
\rep\ \eqref{eq:repr} and Kabluchko \cite{kabluchko:2009} proved this
fact for any max-stable process on $T$ for sufficiently rich measure spaces
$(E,{\mathcal E},\nu)$. In what follows, we will always assume that
the considered max-stable processes have \rep\  \eqref{eq:repr}.
\par
Next we give some basic properties of a stationary max-stable process. 
\bpr\label{prop:1}
The following
statements hold for the skeleton process $(X_t)_{t\in\bbz}$ of the
process \eqref{eq:repr}.
\begin{enumerate}
\item[\rm (1)]
The  \fidi s of $(X_t)$ are 
\regvary\ with index $\alpha$ and the limit \ms s $\mu_h$ of the \fidi s are
given by its values on the complements of the 
rectangles $({\bf0},\bfx]=\{\bfy\in\bbr^h: 0<y_i\le
x_i\,,i=1,\ldots,h\}$, $h\ge 1$, $\bfx=(x_1,\ldots,x_h)$ with $x_i>0$,
$i=1,\ldots,h$:
\beam\label{eq:ooo}
\mu_h\big(
({\bf0},\bfx]^c\big)=
\dfrac{\int_E \max_{i=1,\ldots,h}\big(\dfrac{f_{i}(y)}{x_i}\big)^\alpha\,
\nu(dy)}{\int_E
  f_0^\alpha (y)\,\nu(dy)}\,.
\eeam
\item[\rm (2)]
The \seq\ $(X_t)$ has extremal index $\theta$ \fif\ the limit
\beam\label{eq:lll}
\theta=\lim_{\nto} \dfrac 1 n \dfrac{\int_E\max_{t=1,\ldots,n}f_t^\alpha(y)\,
\nu(dy)}{\int_E
  f_0^\alpha (y)\,\nu(dy)}
\eeam
exists. 
\item[\rm (3)]
The extremogram for the sets $A=(a,\infty)$ and
$B=(b,\infty)$, $a,b>0$, is given by
\beam\label{eq:extrab}
\gamma_{AB}(h)=\dfrac{
\int_E f_0^\alpha(y) \wedge \big(\frac ab
  f_h(y)\big)^\alpha\,\nu(dy)}{a^\alpha\int_E f_0^\alpha(y)\nu(dy)}\,,\quad
h\ge 0\,,
\eeam 
and, for $a=b=1$,
\beam\label{eq:extraba}
\rho(h)= \dfrac{\int_E f_0^\alpha(y) \wedge 
  f_h^\alpha(y)\,\nu(dy)}{
\int_E f_0^\alpha(y)\,\nu(dy)}\,,\quad h\ge 1\,.
\eeam
\item[\rm (4)]
Let $S_1,S_2$ be finite disjoint subsets of $\bbz$
and $\sigma(C)$ the $\sigma$-field generated by 
$(X_t)_{t\in C}$ for any $C\subset \bbz$. 
Recall the
$\alpha$-mixing coefficient  relative to the sets $S_1,S_2$.
\beao
\alpha(S_1,S_2) =\sup_{A\in \sigma(S_1),B\in\sigma(S_2)}|P(A\cap B)- 
P(A)\,P(B)|\,. 
\eeao
and for $S_{-\infty}^0=\{\ldots,-1,0\}$,
$S_h^\infty=\{h,h+1,\ldots\}$, $h\ge 1$, introduce the {\em mixing
  rate \fct }
\beao
\alpha_h=\alpha(S_{-\infty}^0,S_h^\infty)\,,\quad h\ge 1\,.
\eeao
Then there exists a universal constant $c>0$ \st
\beam\label{eq:alpha}
\alpha_h\le c \sum_{s_1=-\infty}^0\sum_{s_2=0}^\infty 
 \int_E 
f_0^\alpha(y)\wedge f_{h+s_2}^\alpha(y)\,\nu(dy)\,,\quad h\ge 1\,.
\eeam
\end{enumerate}
\epr
\bre
If the limit in \eqref{eq:lll} exists it belongs to the interval
$[0,1]$.
Indeed, by stationarity of $(X_t)$,
\beao
\int_E\max_{t=1,\ldots,n}f_t^\alpha(y)\,\nu(dy)\le \sum_{t=1}^n 
\int_E f_t^\alpha(y)\,\nu(dy)= n\,\int_E
f_0^\alpha(y)\,\nu(dy)\,.
\eeao
\ere
\bre\label{rem:2}
Part (4) is a con\seq\ of Corollary 2.2 in Dombry and 
Eyi-Minko \cite{dombry:eyi-minko:2012} proved for $\beta$-mixing.
If  $\int_E f_0^\alpha(y)\wedge f_{h}^\alpha(y)\,\nu(dy)
\le c_0\,\ex^{-c_1\,h}$\,,
$h\ge 1$, for some constants $c_0,c_1>0$, then we conclude that 
$\alpha_h\le C\,\ex^{-c_1\, h}$,  for some $C>0$.

\begin{figure}[htbp]\label{fig:1}
\centerline{
\epsfig{figure=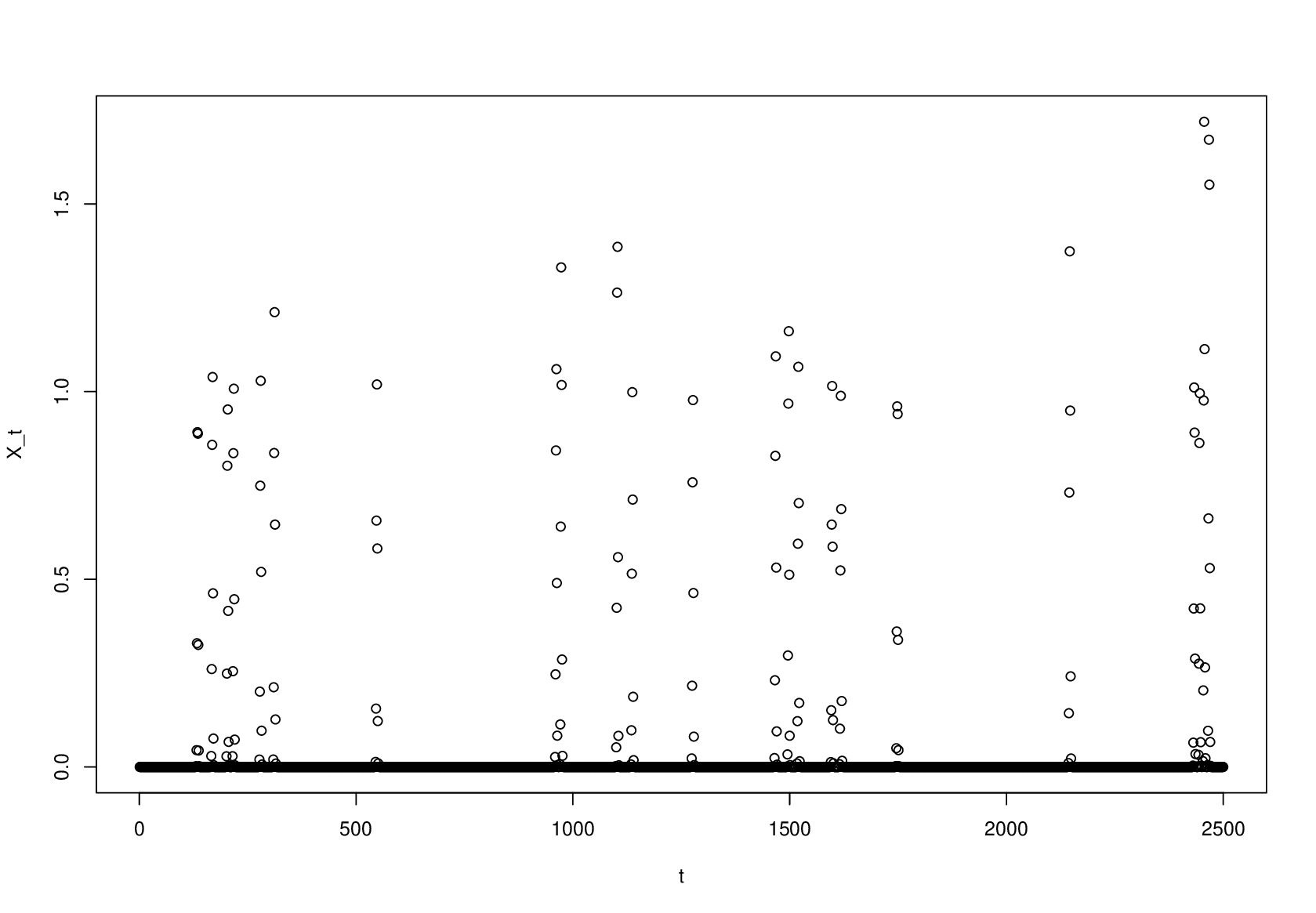,height=7cm,width=7cm,angle=-90}}
\bfi{\em 
Max-stable process \eqref{eq:maxstab} where $f$ is the
standard normal density and $\alpha=5$. Extremal clusters 
are clearly visible. }
\efi
\end{figure}

\ere
\begin{proof} {\bf Part (1)}
Since the integrals $\int_E^{\vee} f_{i} d M_\nu^\alpha$,
$i=1,\ldots,h$,
are supported on $(0,\infty)$ it suffices to show that 
there exists a non-null
Radon \ms\ $\mu_{h}$ on the Borel $\sigma$-field of $\ov \bbr^d_0\cap
(0,\infty]^d$ 
\st\
\beam\label{eq:ttt}
n\,P\big(a_n^{-1} (X_1,\ldots,X_h))\in [{\bf0},\bfx]^c\big)
\to \mu_{h}\big([{\bf0},\bfx]^c\big)\,,
\eeam
where $\bfx$ is chosen \st\ $[{\bf0},\bfx]^c$ is a
$\mu_{h}$-continuity set and 
\beao
P(X>a_n)=1-\exp\big\{-a_n^{-\alpha}\int_E
  f_0^\alpha (x)\nu(dx)\big\}\sim n^{-1}\,,\eeao 
see Resnick \cite{resnick:2007},
Theorem 6.1. 
A Taylor expansion argument
shows that we can always choose 
\beao
a_n= n^{1/\alpha} \big(\int_E
  f_0^\alpha (x)\,\nu(dx)\big)^{1/\alpha}\,.
\eeao
An application of \eqref{eq:fidi} and a Taylor expansion yield  
\eqref{eq:ttt} with limit as specified in \eqref{eq:ooo}.\\
{\bf Part (2)}
Applying \eqref{eq:fidi} for
$x_i=x>0$, we obtain
\beao
P(a_n^{-1}M_n\le x)&=&
\exp\Big\{-a_n^{-\alpha}x^{-\alpha}\int_E \max_{t=1,\ldots,n} f_t^\alpha(y)
\nu(dy)\Big\}\,.
\eeao
By definition of the extremal index, the \rhs\ must converge to 
$\Phi_\alpha^\theta(x)$ for some  
$\theta\in [0,1]$. Equivalently, the limit $\theta$ in \eqref{eq:lll}
exists. \\
{\bf Part (3)}
As regards the extremogram for sets $A=(a,\infty)$, $B=(b,\infty)$,
$a,b>0$, we have the relation
\beao
P(X_h>b\, x, X_0>a\, x)&=& P(X_h> b \,x)+ P(X_0> a\,x)
-P((X_h/b)\vee (X_0/a)> x)\\
&=& 1-\exp\Big\{ -x^{-\alpha}\int_E \big(\dfrac{f_0(y)}{a}\big)^\alpha\,\nu(dy)\Big\}
-\exp\Big\{ -x^{-\alpha}\int_E \big(\dfrac{f_h(y)}{b}\big)^\alpha\,\nu(dy)\Big\}\\
&&+\exp\Big\{-x^{-\alpha}\int_E \big(\dfrac{f_0(y)}{a}\big)^\alpha\vee 
\dfrac{f_h(y)}{b}\big)^\alpha\,\nu(dy)\Big\}\,.
\eeao
In view of stationarity, 
$\int_E f_h^\alpha(y)\,\nu(dy)= \int_E f_0^\alpha(y)\,\nu(dy)$. Using
a Taylor expansion as $\xto$, we obtain the desired formulas
\eqref{eq:extrab} and  \eqref{eq:extraba}\\
{\bf Part (4)} We obtain from Corollary 2.2 in Dombry and 
Eyi-Minko \cite{dombry:eyi-minko:2012}  for any disjoint closed
countable subsets $S_1$, $S_2$ of $\bbr$
\beao
\alpha(S_1,S_2)\le c \sum_{s_1\in S_1}\sum_{s_2\in S_2} \int_E 
f_0^\alpha(y)\wedge f_{|s_1-s_2|}^\alpha(y)\,\nu(dy)\,.
\eeao
Then \eqref{eq:alpha} is immediate.
\end{proof}
Next we consider two popular models of max-stable processes.
\bexam\label{exam:brown:resnick}{\rm The {\em Brown-Resnick process} (see \cite{brown:resnick:1977}) 
has \rep
\beam\label{eq:brown}
X_t= \sup_{i\ge 1}\Gamma_i^{-1/\alpha}\,
\ex^{W_i(t)-0.5\,\sigma^2(t)}\,,
\quad t\in \bbr\,,
\eeam
where $(\Gamma_i)$ is an enumeration of the points of a unit rate homogeneous
Poisson process on $(0,\infty)$  independent of the iid \seq\ $(W_i)$
of sample continuous mean zero Gaussian processes on
$\bbr$ with stationary increments and variance \fct\ $\sigma^2$.
The max-stable process \eqref{eq:brown} is stationary (Theorem 2
in Kabluchko et al. \cite{kabluchko:schlather:haan:2009}) and its 
\ds\ only depends on the variogram $V(h)= \var(W(t+h)-W(t))$,
$t\in\bbr,h\ge 0$. 
It follows from Example 2.1 in Dombry and Eyi-Minko 
\cite{dombry:eyi-minko:2012} that the \fct s $(f_t)$ in \rep\ 
\eqref{eq:repr} satisfy the condition
\beao
\int_E 
f_0^\alpha(y)\wedge f_{h}^\alpha(y)\nu(dy)\le c \,\ov \Phi(0.5
\sqrt{V(h)})\,,
\eeao
where $\Phi$ is the standard normal \ds . For example, if $W$ is
standard \BM , $V(h)=h$, $\ov \Phi(0.5\,\sqrt{h})\sim c\, \ex^{-h/8}
h^{-0.5}$, as $h\to\infty$.
An application of Remark~\ref{rem:2} shows that $(\alpha_h)$ decays at an 
exponential rate. 
\par
Recently, the Brown-Resnick process has attracted some attention
for modeling spatio-temporal extremes; see 
\cite{kabluchko:2009,kabluchko:schlather:haan:2009,stoev:2008,oesting:kabluchko:schlather:2012}.
The processes \eqref{eq:brown} can be extended to
  random fields on $\bbr^d$. These fields found various applications for
  modeling spatio-temporal extremal effects; see Kabluchko et al. 
\cite{kabluchko:schlather:haan:2009}. For further
  spatio-temporal 
applications
  of
max-stable random fields, see also Davis et al. \cite{davis:kluppelberg:steinkohl:2012}.

}
\eexam
\bexam\label{exam:dehaan:pereira}{\rm We consider de Haan and Pereira's \cite{haan:pereira:2006}
 {\em max-moving process}
\beam\label{eq:maxstab}
X_t=\sup_{i\ge 1} \Gamma_i^{-1/\alpha} f(t-U_i)\,,\quad t\in\bbr\,,
\eeam
where $f$ is a continuous Lebesgue density on $\bbr$ \st\ $\int_\bbr
\sup_{|h|\le 1} f(x+h)\,dx<\infty$ and 
$\sum_{i=1}^\infty\vep_{(\Gamma_i,U_i)}$ are the points of a unit rate
homogeneous 
Poisson random \ms\ on $(0,\infty)\times \bbr$. 
\par
The resulting process $(X_t)$ is $\alpha$-max-stable
and stationary. According to Example 2.2 in  Dombry and Eyi-Minko 
\cite{dombry:eyi-minko:2012},
\beao
\int_E 
f_0^\alpha(y)\wedge f_{h}^\alpha(y)\nu(dy)\le c \int_\bbr
\min(f(-x),f(h-x))\,dx\,,\quad h\ge 0\,.
\eeao
For example, if $f$ is the standard normal density, this implies
that $(\alpha_h)$ decays to zero faster than exponentially, i.e. the
memory in this \seq\ is very short.
In Figure~\ref{fig:1} a simulation of the corresponding process 
\eqref{eq:maxstab} for $\alpha=5$ is shown.
}
\eexam
\section{Estimation of the extremogram}\setcounter{equation}{0}
\subsection{Asymptotic theory}\label{subsec:21}
Natural estimators of the extremogram are obtained by replacing the
\pro ies in the limit relations \eqref{eq:extrem}  and 
\eqref{eq:rho} with their empirical counterparts. 
In this context, one works with 
quantities which are derived from the {\em tail empirical process};
see the monographs de Haan and Ferreira \cite{dehaan:ferreira:2006}, 
Resnick \cite{resnick:2007} for the underlying theory. For the
introduction of the sample extremogram, consider an $\bbr^d$-valued
strictly stationary \regvary\ process $(X_t)$ and a Borel set $C\subset \ov
\bbr_0^d$ bounded away from zero. Then, for any \seq\ $m=m_n\to\infty$
with $m_n/n\to 0$ as $\nto$, we define the following estimator
of $P_m(C)=m\,P(a_m^{-1} X\in C)$:
\beao
\wh P_m(C)=\dfrac m n \sum_{t=1}^n I_{\{a_m^{-1} X_t\in C \}}
\eeao
A possible choice of $(a_m)$ is given by $P(|X|>a_m)\sim
m^{-1}$.
By definition of \regvar\ of $X$, for any $\mu_1$-continuity set
$C$, 
\beao
E[\wh P_m(C)]= m\, P_m(C)\to \mu_1(C)\,.
\eeao
Here the condition $m_n\to\infty$ as $\nto$ was crucial for \asy\ unbiasedness. For the
calculation of the \asy\ variance of $\wh P_m(C)$ we   
assume  the following condition: 
\begin{enumerate}
\item[(M)]
The \seq\ $(X_t)$ is $\alpha$-mixing with rate \fct\ $(\alpha_h)$
and there exists a \seq\ $r_n\to\infty$ \st\ $r_n/m_n\to 0$ as $\nto$,
\beam\label{eq:mix1}
\lim_{\nto} m_n \sum_{h=r_n}^\infty \alpha_h=0
\eeam
and for every $\epsilon>0$, 
\beam\label{eq:mix2}
\lim_{\kto} \limsup_{\nto} m_n\sum_{h=k}^{r_n} P(|X_h|>\epsilon\,
a_m\,,|X_0|>\epsilon\,a_m)=0\,.
\eeam
\end{enumerate}
This condition is technical: \eqref{eq:mix1}
imposes some rate on the mixing \fct\ $(\alpha_h)$ and \eqref{eq:mix2}
avoids ``extremal long range dependence'';  \eqref{eq:mix2} is an
\asy\
independence condition in the spirit of the classical 
condition $D'$; see Leadbetter et al. 
\cite{leadbetter:lindgren:rootzen:1983}, 
Embrechts et al. \cite{embrechts:kluppelberg:mikosch:1997}.
The quantities $m_n$ and $r_n$ have some straightforward
interpretation as size in large-small block scheme: 
the sample $X_1,\ldots,X_n$ consists of roughly $[n/m_n]$ 
large disjoint blocks of size $m_n$. After chopping off the first
$r_n$ elements in each large block one aims at  ensuring the
\asy\ independence of the resulting large blocks.   
\par
If $C$ is a $\mu_1$-continuity set and $C\times \ov
\bbr_0^{d(h-1)}\times C$ are $\mu_{h+1}$-continuity sets for every
$h\ge 1$, \regvar\ of $X$ implies
\beam\label{eq:asymp}
\var[\wh P_m(C)]\sim \dfrac{m}{n} \,V(C)\,,
\eeam
where
\beao%\label{eq:V}
V(C)=\mu_1(C)+ 
2\,\sum_{h=1}^\infty \tau_h(C)\quad\mbox{and}\quad \tau_{h}(C)= 
\mu_{h+1}( C\times \ov \bbr_0^{d(h-1)}\times C)\,,\quad h\ge 1\,,
\eeao 
and we also assume that the infinite series is finite.
The \asy\ relation \eqref{eq:asymp} indicates that the condition 
$m_n/n\to 0$ is needed to ensure the consistency of the estimator 
$\wh P_m(C)$. Under additional conditions, $(\wh P_m(C))$ is \asy ally
normal and this property also holds jointly for finitely many
sets $C_1,\ldots,C_h$. The complicated form of the \asy\ variance
in \eqref{eq:asymp} suggests that it is difficult to apply this \clt\
for constructing \asy\ confidence bands.
\par
The motivating examples of extremograms are limits of
conditional \pro ies
\beao
\rho_{AB}(h) = \lim_{\xto} P(x^{-1} X_h\in B\mid x^{-1} X_0\in
A)\,,\quad h\ge 0\,, 
\eeao
for sets $A,B$ bounded away from zero. In
Section~\ref{subsec:extremogram}
we mentioned the close relationship of
$\rho_{AB}$ with a cross-correlation \fct .
Replacing the \pro ies in these
conditional \pro ies by estimators of the type $\wh P_m(C)$ and
applying the corresponding central limit theory from
\cite{davis:mikosch:2009c}, Section 3, and the \cmt , one obtains an \asy\ theory for
the {\em ratio estimators}
\beao
\wh \rho_{AB}(h)=
\dfrac{\sum_{t=1}^{n-h} I_{\{a_m^{-1}X_t\in A\,,a_m^{-1}X_{t+h}\in B
    \}}}
{\sum_{t=1}^nI_{\{a_m^{-1} X_t\in A\}}}\,,\quad h\ge 0\,.
\eeao
The latter estimators only depend on the high threshold
$a_m$ which one typically chooses as a high empirical quantile of the  
data. These estimators can be interpreted as a sample cross-correlation
\fct .
\par
We recall a \clt\ for these estimators; see Corollary 3.4 in
\cite{davis:mikosch:2009c} and its correction Theorem 4.3 
in \cite{davis:mikosch:2012a}.
\bth\label{thm:1}
Let $(X_t)$ be an $\bbr^d$-valued  strictly stationary \regvary\
\seq\ with index $\alpha>0$.
Assume that the following conditions are satisfied.
\begin{itemize}
\item
The Borel sets $A,B\subset \ov \bbr_0^d$ are bounded away from zero and $\mu_1(A)>0$.
\item
The sets $A$, $B$ are continuous \wrt\ $\mu_1$. 
\item 
Condition {\rm (M)}, $(n/m_n)\alpha_{r_n}\to 0$ as $\nto$.
\item
$m_n=o(n^{1/3})$ or 
\beam\label{eq:no2}
\dfrac{m_n^4} n \sum_{j=r_n}^{m_n} \alpha_j \to 0\quad\mbox{and}
\quad \dfrac{m_n r_n^3}{n}\to 0\quad \mbox{as $\nto$.}
\eeam
\end{itemize}
Then the following \clt\ holds for $h\ge 0$
\beam\label{eq:preas}
\sqrt {\dfrac n {m_n}}\,
\Big[\wh \rho_{AB}(h)- \rho_{AB,m}(h)\Big]_{h=0,\ldots,m}
\std N({\bf0}, (\mu_1(A))^{-4} {\bf \Sigma})\,.
\eeam
for some matrix $\bf\Sigma$,\footnote{This matrix  is complicated 
and
irrelevant for our purposes; see \cite{davis:mikosch:2009c}, (3.15) and
(3.16) for its value.} where $\rho_{AB,m}(h)=P(a_m^{-1} X_h\in B\mid
a_m^{-1} X_0\in A)$.
\ethe
Some comments are here in place.
\begin{itemize}
\item 
The conditions of this result are rather technical but the mixing 
and anti-clustering conditions can be verified 
for standard time series models. For example, if $(X_t)$ is
$\alpha$-mixing with geometric rate, then one can simply choose
\seq s $r_n=[C\log n]$ for a large constant $C>0$ or $r_n=n^\epsilon$,
and $m_n=n^{2\epsilon}$ for suitable small $\epsilon>0$.
\item
The \asy\ variance is not of practical use. Therefore Davis et al. 
\cite{davis:mikosch:2012} suggest an alternative way of constructing
confidence bands for $\wh \rho_{AB}(h)$, by using the stationary
bootstrap introduced by Politis and Romano \cite{politis:romano:1994}. 
\item
Central limit theory and bootstrap consistency for the sample
extremogram do not
follow from standard results for \seq s of mixing stationary \seq s.
This is due to the fact that we deal with 
\seq s of indicator \fct s $(I_{\{a_m^{-1} Y_t\in C\}})$ for certain 
strictly stationary \seq s $(Y_t)$ and sets $C$ bounded away from
zero. The \seq s $(I_{\{a_m^{-1} Y_t\in C\}})$ constitute a {\em triangular array of row-wise strictly
stationary \seq s} for which, to the best of our 
knowledge, standard \asy\ theory 
is not available. 
\item
The \con\ rate $\sqrt{n/m}$ 
in the \clt\ \eqref{eq:preas} is due to the triangular array
structure;
it can be significantly slower than standard $\sqrt{n}$-rates. 
\item
We call $\rho_{AB,m}$ in \eqref{eq:preas} the {\em pre-\asy\ extremogram } since, in
general, one cannot replace the centering constants $\rho_{AB,m}(h)$
by their limits $\rho_{AB}(h)$; see Example~\ref{exam:1} below. 
Moreover, it is in general very
difficult to show that 
\beam\label{eq:exam1}
\sqrt{\dfrac{n}{m_n}} |\rho_{AB,m}(h)-\rho_{AB}(h)|\to
0\,,\quad \mbox{as $\nto$,}
\eeam
 even for ``nice'' models such as the \garch . For
this well studied model one lacks precise information about the tail
behavior. The \clt\ \eqref{eq:preas} (and
related bootstrap procedures) are then used to approximate  the
conditional \pro ies $\rho_{AB,m}(h)$. These have a very concrete
interpretation in contrast to their less intuitive limits
$\rho_{AB}(h)$.
\item
If \eqref{eq:exam1} holds with a rather slow \con\ rate one
faces a bias problem. This problem can be observed e.g. 
for a simulated \sv\ process $(X_t)$ and $A=B=(1,\infty)$. Then
$\rho(h)=\rho_{AA}(h)=0$ for $h\ge 1$. If $(X_t)$ is 
$\alpha$-mixing with geometric rate
it can be verified that \eqref{eq:exam1} and \eqref{eq:no2} hold for
$m_n=n^\gamma$, $\gamma\in (1/3,1)$, and then  \eqref{eq:preas}
applies with
$\rho_{AA,m}(h)$ replaced by $\rho(h)=0$; see
\cite{davis:mikosch:2009c}, Section 4.2.  Also notice that $\wh
\rho_{AA}(h)$ is of the order $1/m$ in the iid case and hence greater
than zero.
\item 
The formulation of the results in 
\cite{davis:mikosch:2009c,davis:mikosch:2012,davis:mikosch:2012a}
related to Theorem~\ref {thm:1} involve 
various other continuity conditions on sets.  These conditions
can be avoided as a close inspection of
the proofs shows: these conditions follow from continuity of 
$A$ and $B$ \wrt\ $\mu_1$. Indeed, one needs that sets of the 
form $\bigotimes_{i=1}^k C_i$ are $\mu_k$-continuity sets, where
$C_i\in \{A,B,\ov \bbr_0^d,A\cap B\}$, $i=1,\ldots,k$, and at least
one of the sets $C_i$ does not coincide with $\ov\bbr_0^d$. Let 
$S$ be the set of indices $i$ \st\ $C_i$ does not coincide with
$\ov\bbr_0^d$. Then
\beao
\partial\Big(\bigotimes _{i=1}^k C_i\Big)\subset \bigcup_{i\in S} 
(\ov \bbr_0^d)^{i-1} \times\partial  C_i\times 
(\ov \bbr_0^d)^{k-i-1}\,. 
\eeao
The sets in the union have $\mu_k$-\ms\ zero. For the sake of
argument assume that $i=1\in S$. Then
\beao
\lim_{\nto}n\,P(a_n^{-1} (X_1,\ldots,X_k)\in \partial C_1\times
(\ov\bbr^d_0)^{k-1})&=&\lim_{\nto}n\,P(a_n^{-1} X_1\in \partial
C_1)\\&=&\mu_1(\partial C_1)=0\,.
\eeao
\item
In applications one needs to choose the threshold $a_m$ in some
reasonable way. The choice of a threshold is inherent to extreme value
statistics to which no easy solution exists. In
\cite{davis:mikosch:2009c} we advocated to choose $a_m$ as a fixed high/low
empirical quantile of the absolute values of the data and to
experiment with several quantile values. If the plot of the sample
extremogram is robust for a range of such quantiles one can choose 
a quantile from that region. The theory in
\cite{davis:mikosch:2009c,davis:mikosch:2012,mikosch:zhao:2012}
is based on {\em deterministic} values $a_m$. The heuristic method
described above advocates the choice of a {\em data dependent}
threshold.  Recent work by Kulik and
Soulier  \cite{kulik:soulier:2011} yields  a theory for a modified 
sample extremogram
with data dependent threshold in the case of  short and long memory \sv\ 
processes.   
\end{itemize}
\bexam\label{exam:1} {\em 
In \cite{davis:mikosch:2009c} we did not provide a concrete
example of a \seq\ $(X_t)$ for which \eqref{eq:exam1} does not
hold. Such counterexamples can easily be constructed from 
max-stable strictly stationary processes 
with Fr\'echet marginals; see  Section \ref{subsub:examples}. We
assume the conditions of Proposition~\ref{prop:1}.
For simplicity choose $\alpha=1$, 
$A=B=(1,\infty)$ and $a_m= m \int_E f_0(x)\,\nu(dx)$. The \fct\ $\rho$ is 
given in \eqref{eq:extraba}.
By Taylor expansion, 
\beao
\rho_{AA,m}(h)&=& \dfrac{P(\min(X_0,X_h)>a_m)}{P(X>a_m)}
= \dfrac{1-\ex^{- a_m^{-1} \int_E f_0(x) \wedge
    f_h(x)\,\nu(dx)}}{1-\ex^{-a_m^{-1} \int_E f_0(x)\,\nu(dx)} }\\
&=& \rho(h) - m^{-1} c_h\,(1+o(1))\,, 
\eeao
for some constant $c_h\ne 0$.
Hence 
\beao
(n/m_n)^{0.5}|\rho_{AA,m}(h)-\rho(h)|\sim c_h\, (n/m^3)^{0.5}\,,
\eeao
and the \rhs\ converges to zero \fif\ $n^{1/3}=o(m)$
and the rate of \con\ to zero can be arbitrarily small. The latter
condition is, of course, in contradiction with $m=o(n^{1/3})$
which is one possible sufficient condition for \eqref{eq:preas}.
Fortunately, the other sufficient condition \eqref{eq:no2} can still
be satisfied if $m=o(n^{1/3})$ does not hold. For example, if
$\alpha_h$ decays to zero at a geometric rate  and one chooses
$r_n=[C \log n]$ for some $C>0$ and $m_n=n^\gamma$ for some $\gamma\in
(1/3,1)$. Particular cases with geometric decay of
$(\alpha_h)$ were mentioned
in Examples~\ref{exam:brown:resnick} and 
\ref{exam:dehaan:pereira}.
}
\eexam

  \subsection{Cross-extremogram for bivariate time series}\label{sec:bivext}
While the definition of the extremogram covers the case of 
multivariate time series, it is of limited value if the index of
regular variation is not the same across the component series.  For
example, consider two regularly varying univariate 
strictly stationary time 
series $(X_t)$ and $(Y_t)$ with tail indices $\alpha_X<\alpha_Y$.
Then, assuming $((X_t,Y_t))_{t\in\bbz}$ 
stationary, this bivariate time series 
would be regularly varying with index $\alpha_X$, and
for Borel sets $A,B$ bounded away from zero,  $\wt A= A\times \bbr$ and $\wt B=\bbr\times B$,
$$
\rho_{\wt A\wt B}(h)=\lim_{x\to\infty}P(Y_{h}\in x B \mid X_0\in x
A)=
\lim_{\xto}
P((X_h,Y_h)\in x \wt B\mid (X_0,Y_0)\in x\wt A)
=0\,,\quad h\in\bbz\,.
$$
The \asy\ theory of Section~\ref{subsec:21} is applicable to the
sets $\wt A= A\times \bbr$ and $\wt B=\bbr\times B$.
In this case, no extremal 
dependence between the two series would be \ms d.  
To avoid these rather uninteresting cases and 
obtain a more meaningful measure of extremal dependence, 
we transform the two series so that they have the same marginals.  
In extreme value theory, the transformation to the unit 
Fr\'echet distribution is standard.  
For the sake of argument, assume 
that both $X_t$ and $Y_t$ are positive so that the 
focus of attention will be on extremal dependence in 
the upper tails.  The case of extremal dependence in 
the lower tails or upper and lower tails is similar.  
Under the positivity constraint, if $F_1$ and $F_2$ 
denote the \df s of $X_t$ and $Y_t$, 
respectively, and are continuous, then the two transformed series, 
$\wt X_t=G_1(X_t)$  and $\wt Y_t=G_2(Y_t)$
with  $G_i(z)=-1/\log(F_i(z))$, $i=1,2,$
have unit Fr\'echet marginals $\Phi_1$; see \eqref{eq:frechet}. 
Now assuming that the bivariate time series $((\wt X_t,\wt Y_t))_{t\in\bbz}$ is regularly varying, we define the {\it cross-extremogram} by
\beao
\rho_{\wt A\wt B}(h)=\lim_{x\to\infty}P(\wt Y_{h}\in x\,B\mid \wt X_0 \in
x\,A)\,,\quad
h\in \bbz\,,
\eeao
At first glance, this may seem inconvenient 
since transformations to unit Fr\'echet marginals are required.   
If one restricts attention to sets $A$ and $B$ that are 
intervals bounded away from 0 or finite unions of such sets the
transformation simplifies: if $a_n$ denotes the 
$(1-n^{-1})$-quantile of $\Phi_1$, then  
by monotonicity of $G_i$, $\{\wt X_h\in a_n\,A\}=\{X_h\in
a_{X,n}\,A\}$ and $\{\wt Y_h\in a_n\,B\}=\{Y_h\in a_{Y,n}\,B\}$, 
where $a_{X,n}$ and $a_{Y,n}$ are the respective 
$(1-n^{-1})$-quantiles of the distributions of $X_t$ and $Y_t$.  
For sets $A$ and $B$ of the required form, the cross-extremogram becomes
\beam\label{eq:crossext}
\rho_{\wt A\wt B}(h)&=&\lim_{n\to\infty}P(Y_h\in a_{Y,n}B\mid  X_0 \in a_{X,n}\,A)\,.
\eeam
Thus we do not actually need to find the transformations 
converting the data to unit Fr\'echet, only the component-wise quantiles, $a_{X,n}$ and $a_{Y,n}$, need to be calculated.
Clearly, this notion of extremogram extends to more than two time
series.

\section{An example: Equity indices}\label{sec:EQIndex}
\begin{figure}[ht]
\begin{center}
\centerline{\includegraphics[height=5cm,width=16cm]{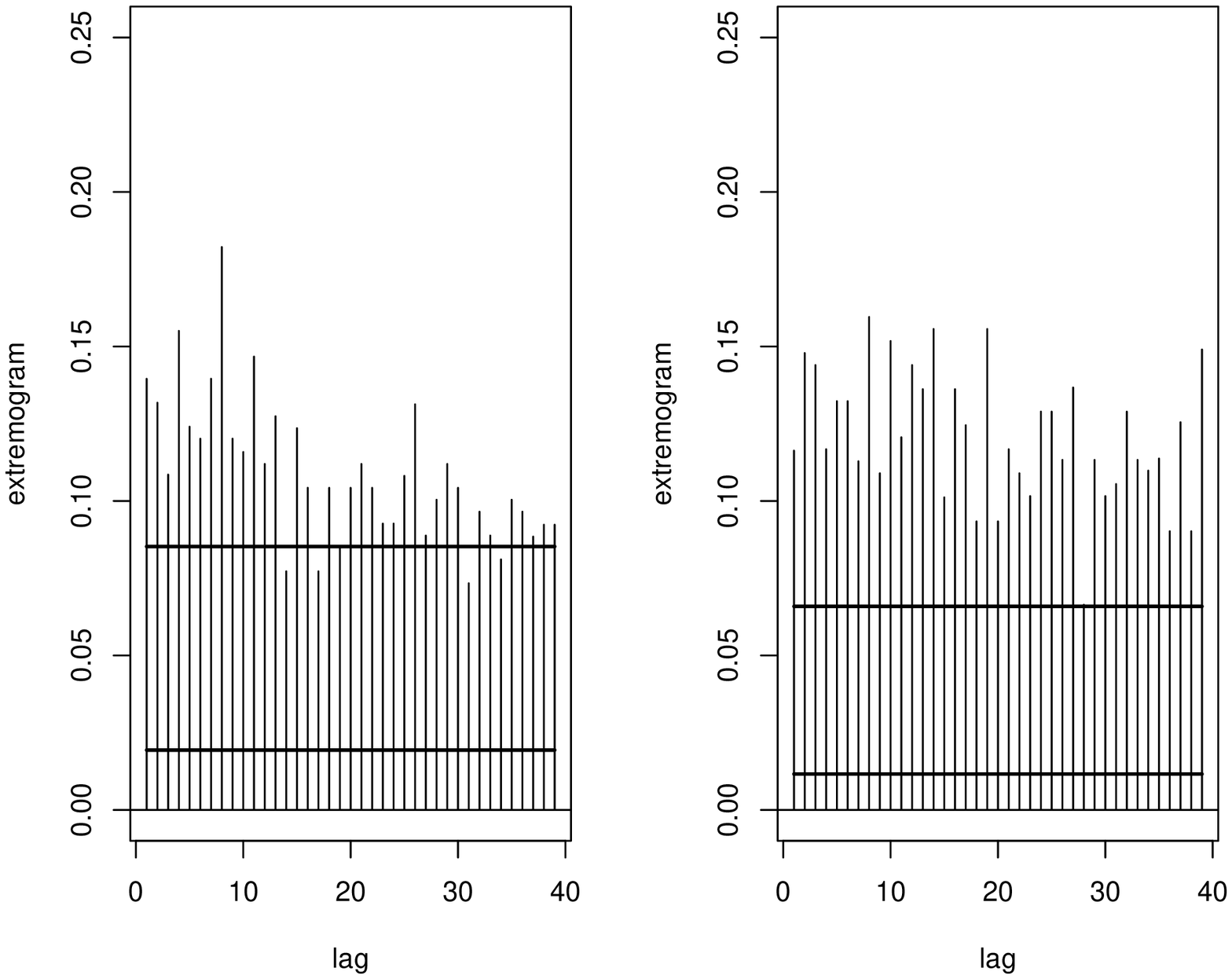}}
\centerline{\includegraphics[height=5cm,width=16cm]{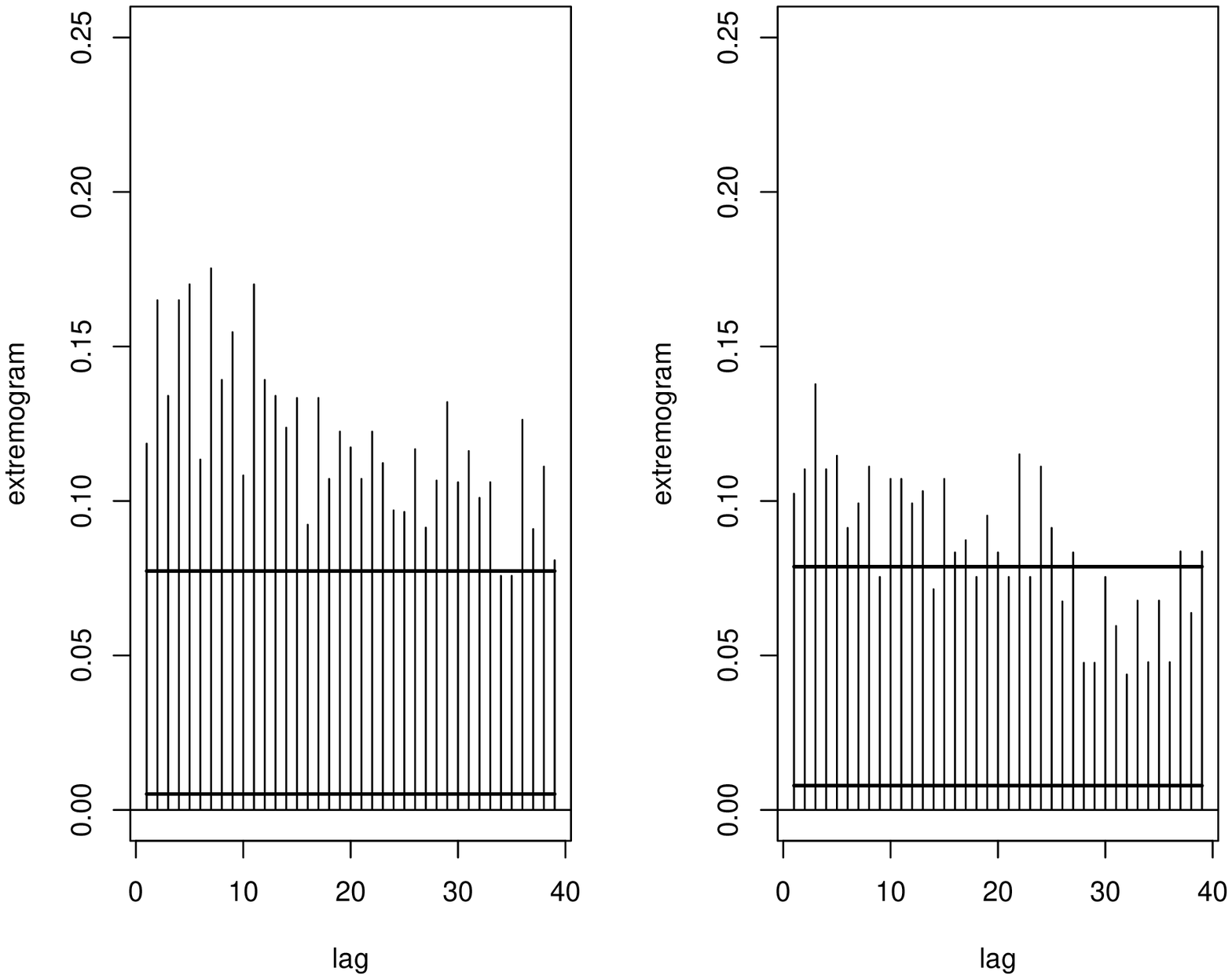}}
\end{center}
\bfi{The sample extremogram for the lower tails of the FTSE {\rm (top left)}, S\&P
  {\rm (top right)}, DAX {\rm (bottom left)} and Nikkei. The solid
  lines represent 98\% confidence bands 
based on $99$ random permutations of the data.}\label{fig:FTSESP}
\efi
\end{figure}

We consider daily log-returns equity indices of 
four countries:  S\&P 500 for US, FTSE 100 for UK, 
DAX for Germany, Nikkei 225 for
Japan. 
Figure ~\ref{fig:FTSESP}
shows the sample extremogram $\wh\rho=\wh\rho_{AA}$ for the negative 
tails ($A=B=(-\infty,-1)$ with $a_m$ estimated 
as the 96\% empirical quantiles of the absolute values of the 
negative data) 
applied to 6,443 
log-returns of the FTSE and S\&P (April 4,
1984
to October 2, 2009), to 4,848
log-returns of the DAX (November 13, 1990
to October 2, 2009) and to
6,333 log-returns of the Nikkei (August 23, 1984 to 
October 2, 2009).\footnote{As noted in the literature, 
the lower tails of returns tend to be heavier 
than the upper tails.  Similar plots (not shown here) of the sample 
extremogram for the upper tails also reveal extremal 
dependence, but to a lesser extent than seen in the lower tails.}
The solid horizontal lines in the plots represent
98\% confidence bands. They correspond to the maximum and minimum 
of the sample extremogram at lag 1 based on 
99 random permutations of the data. If the
data were independent random permutations would not change the
dependence structure: values $\wh \rho_{AB}(h)$ which are
outside the confidence bands indicate that there is significant 
extremal dependence at lag $h$. 
The sample extremograms for all four indices decay 
rather slowly to zero,
with S\&P the slowest.
Among the four indices, the Nikkei displays the least amount 
of extremal dependence as measured by the extremogram.
The top graphs in
Figure~\ref{fig:FTSESP} indicate extremal dependence in the lower 
tail over a period of
40 days.
\par
\begin{figure}[ht]
\begin{center}
\centerline{\includegraphics[height=5cm,width=16cm]{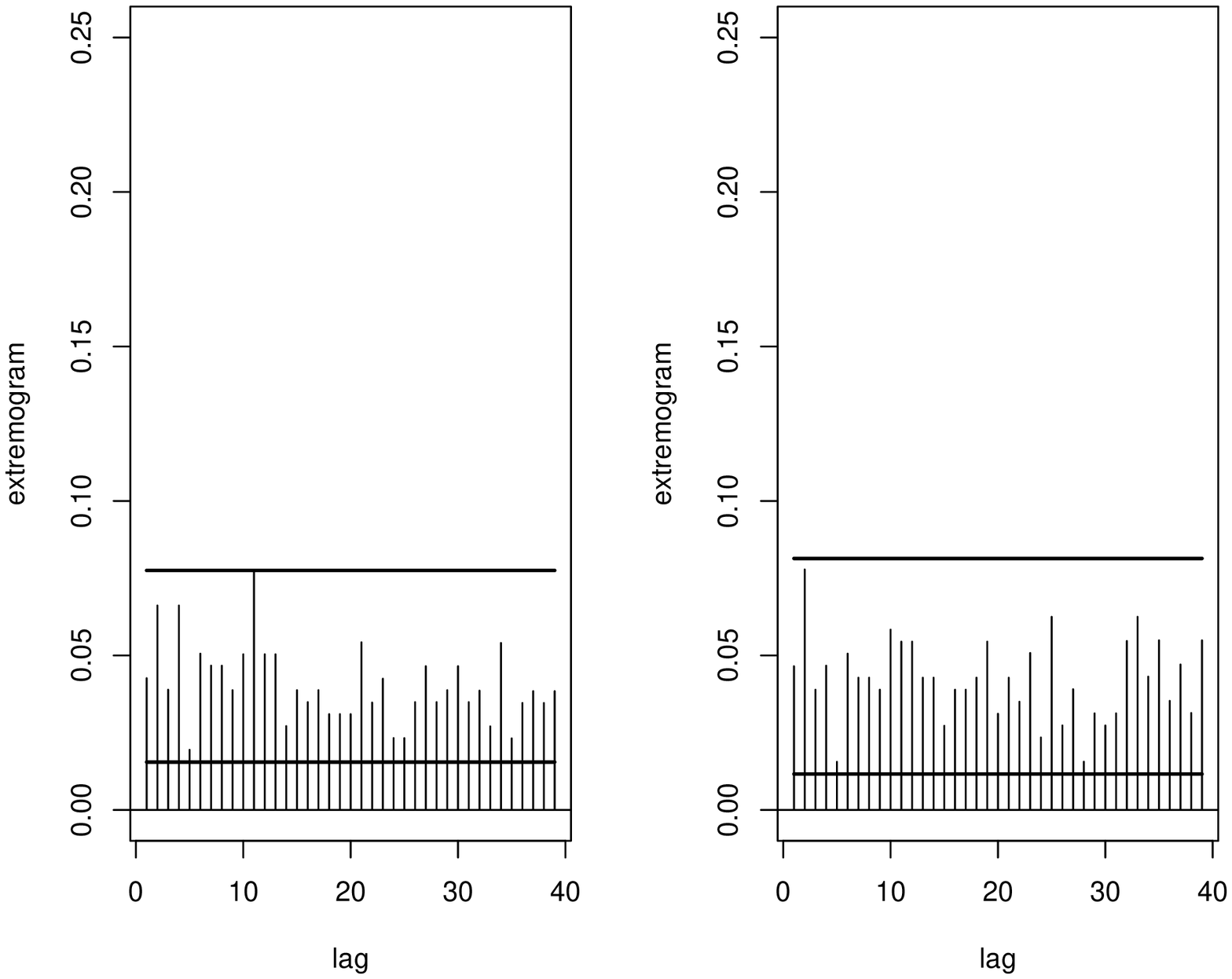}}
\end{center}
\bfi{\small The sample extremograms for the
filtered FTSE {\rm (left)} and filtered S\&P {\rm (right)} series. The
bold lines represent 98\% confidence bounds 
based on $99$ random permutations of the series.}
\label{fig:FTSESPres}
\efi
\end{figure}

 
We assume that the log-return series are modeled by a \garch\ process
(see Section~\ref{subsub:examples} for its definition). Then we 
can estimate its parameters, calculate the volatility \seq\ 
$(\wh \sigma_t)$ and the {\em filtered} \seq\ 
$\wh Z_t=X_t/\wh\sigma_t$, $t=1,\ldots,n$.
Figure~\ref{fig:FTSESPres}
shows the sample extremograms $\wh\rho$ for the  filtered 
FTSE and S\&P \seq s. These plots confirm that much of the 
extremal dependence (as measured by the extremogram) has been removed.
Hence the extremal dependence in the log-returns is due to the
volatility \seq s $(\sigma_t)$, as suggested by the GARCH model.
\par
For a bivariate time series $((X_t,Y_t))_{t\in\bbz}$ the sample
cross-extremogram is given by
\beao
\wh\rho_{\wt A\wt B}(h)=\dfrac{\sum_{t=1}^{n-h}I_{\{Y_{t+h}\in a_{Y,m}\,B,X_t\in a_{X,m}A\}}}{\sum_{t=1}^{n}I_{\{X_t\in a_{X,m}A\}}}\,,
\eeao
where $a_{X,m}$ and $ a_{Y,m}$ are the 
$(1-m^{-1})$-quantiles of the marginal
\ds s of $X$ and $Y$, respectively. For applications, they need to
be replaced by the corresponding
empirical quantiles.
\par
We calculate the sample cross-extremograms for the pairs of
the log-returns series, again for the negative tails, i.e.
$A = B = (-\infty,-1)$, and $a_{X,m}$ and $a_{Y,m}$ are chosen as
the 96\%  empirical quantiles  
of the negative values of the corresponding component samples.
Since the samples have different sizes
we consider those periods for which we have observations on both
indices.  
\par
\begin{figure}[ht]
\begin{center}
\centerline{\includegraphics[height=3.5cm,width=16cm]{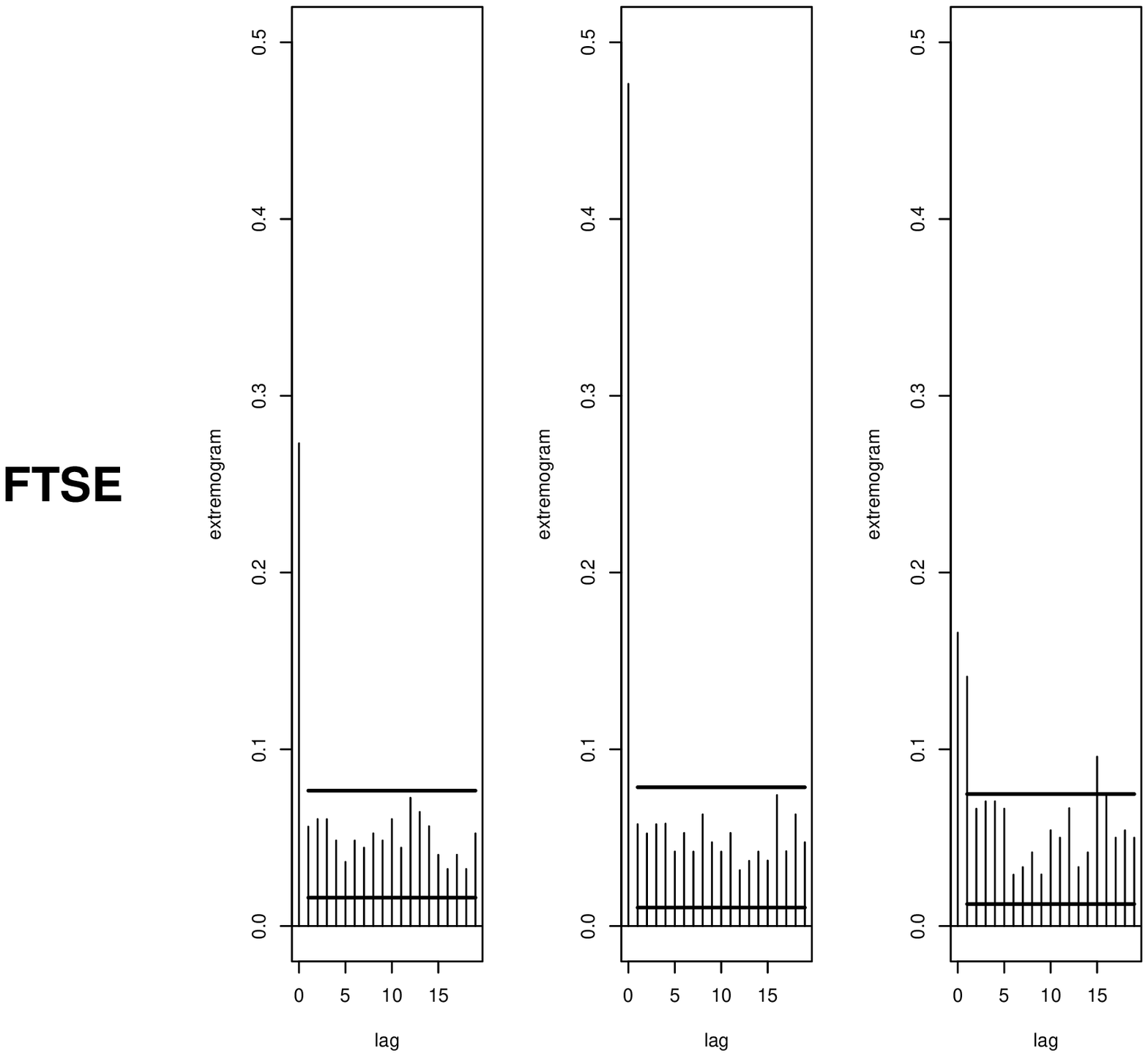}}
\end{center}
\begin{center}
\centerline{\includegraphics[height=3.5cm,width=16cm]{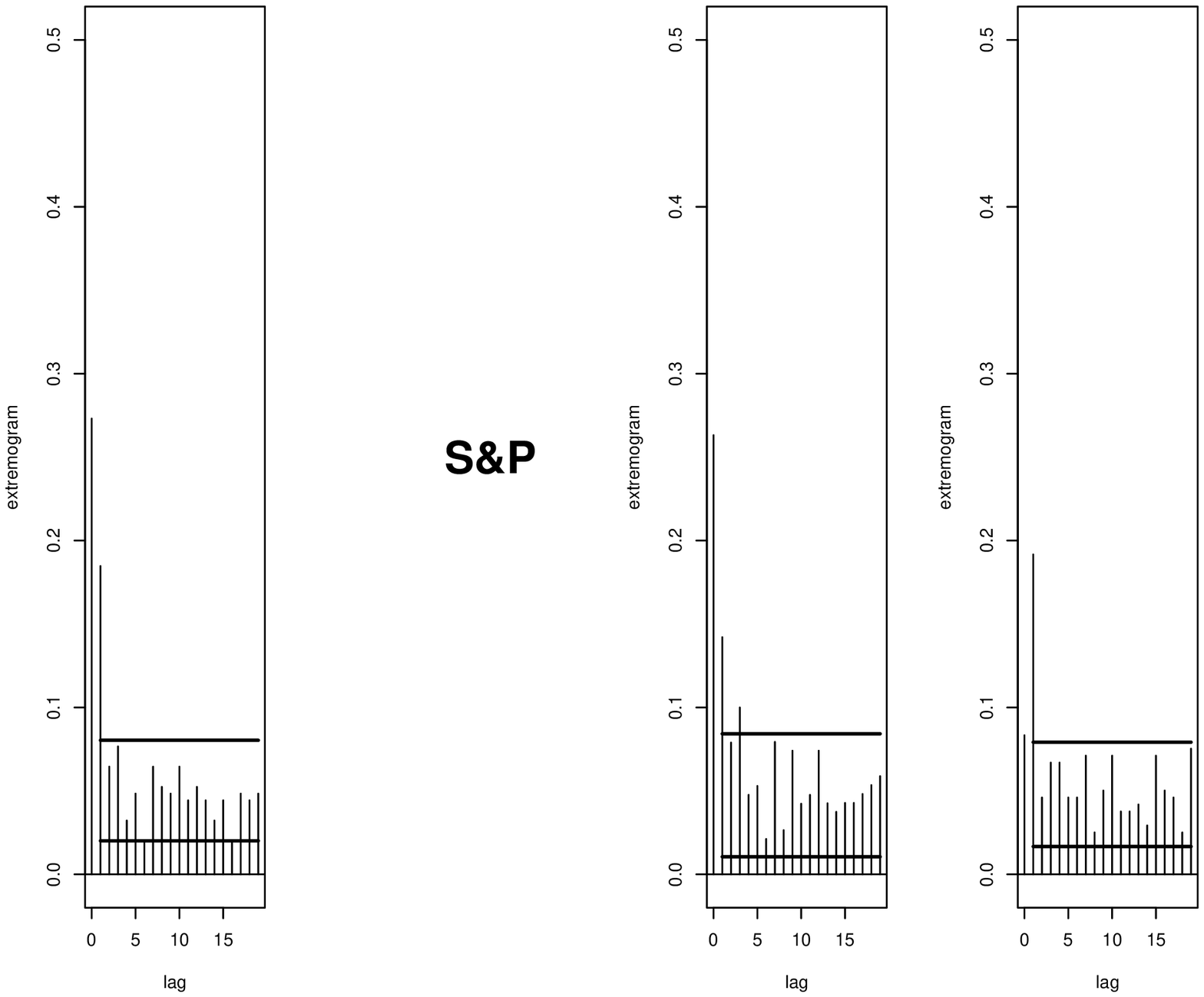}}
\end{center}
\begin{center}
\centerline{\includegraphics[height=3.5cm,width=16cm]{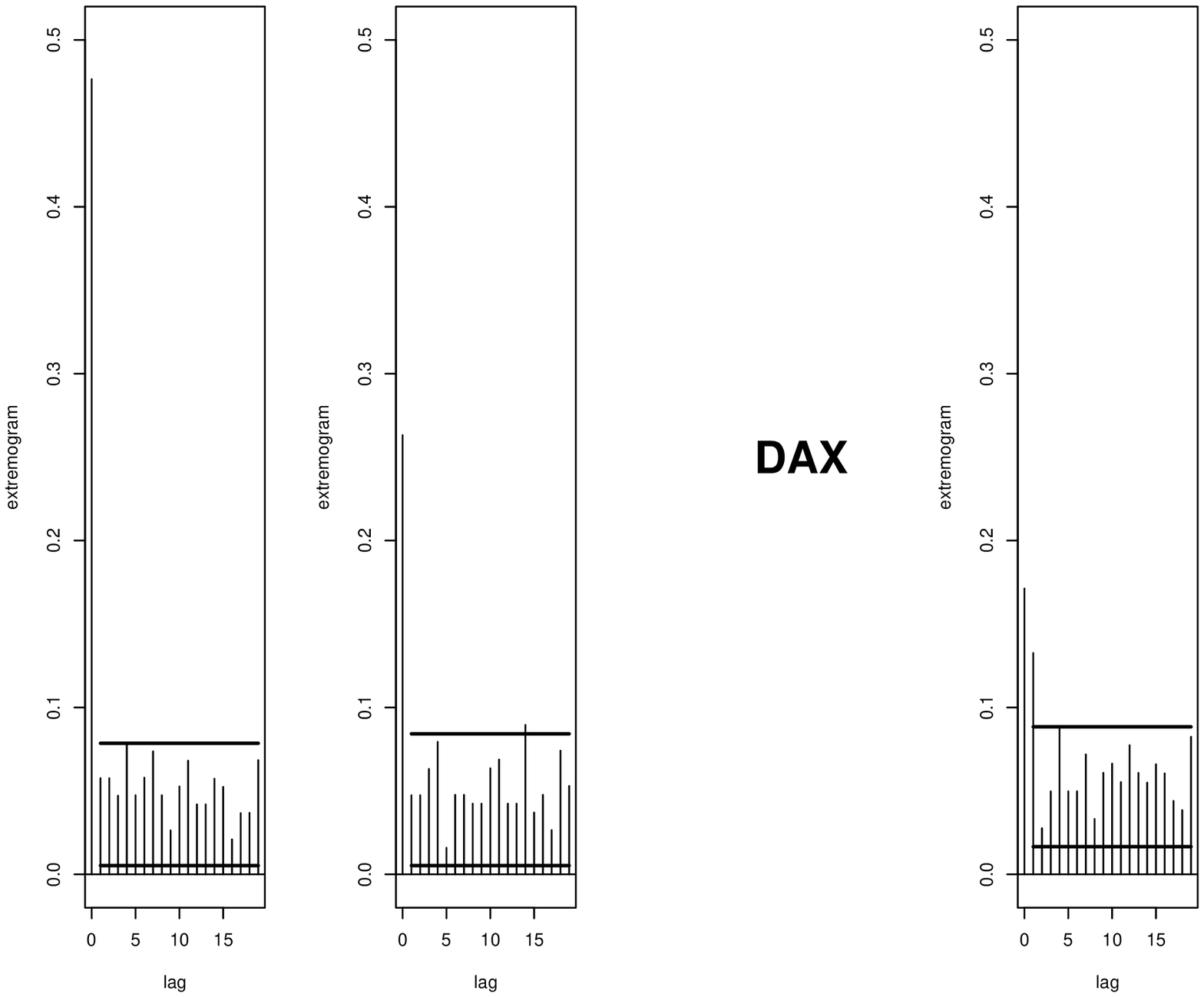}}
\end{center}
\begin{center}
\centerline{\includegraphics[height=3.5cm,width=16cm]{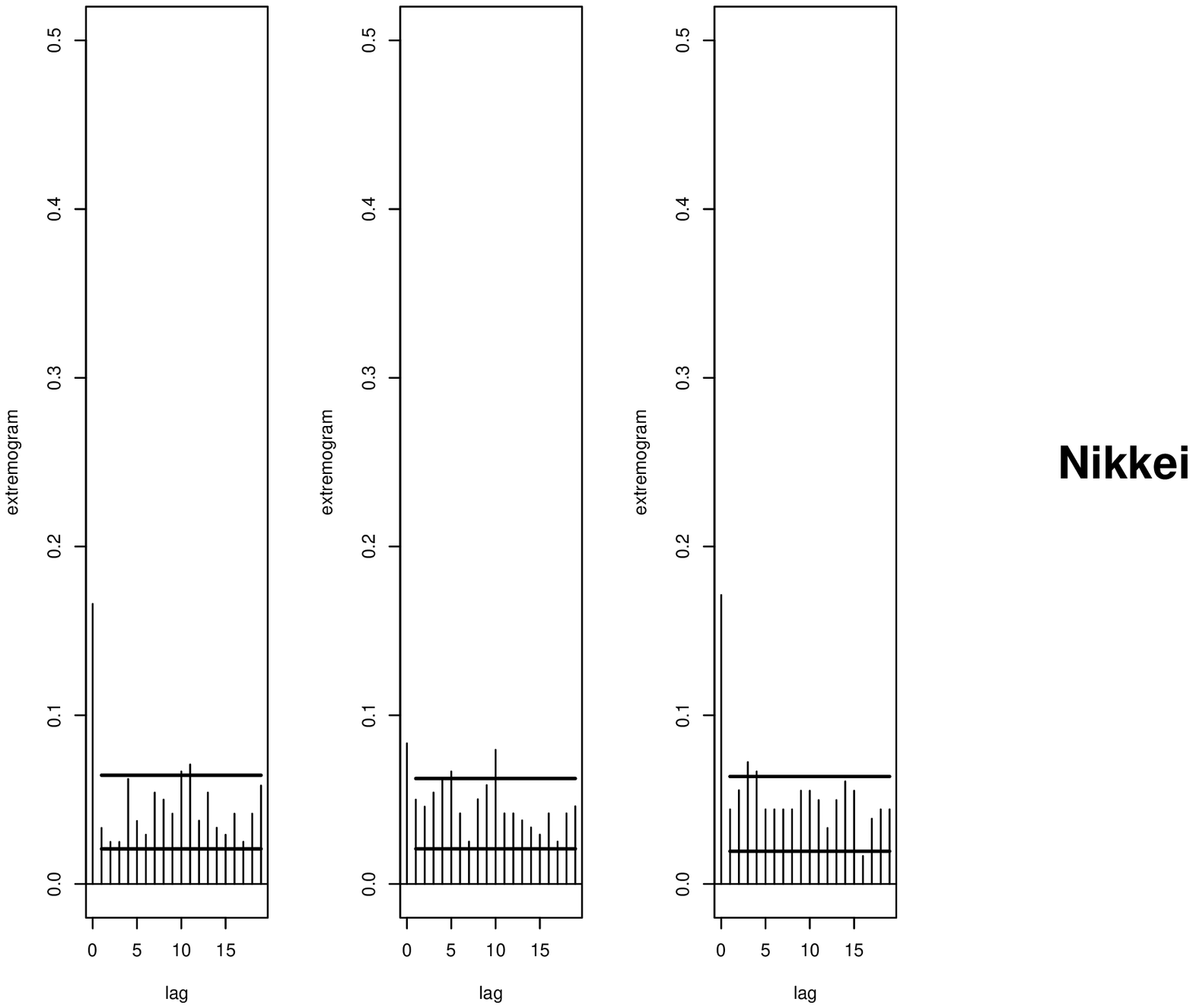}}
\end{center}
\bfi{\small The sample cross-extremograms for the filtered
FTSE, S\&P, DAX and Nikkei series. For the first row, $(X_t)$
is the filtered FTSE and $(Y_t)$ are the filtered S\&P, DAX and Nikkei
{(from left to right)}. For the second, third and fourth rows, the $X_t$'s
are the filtered S\&P, DAX and Nikkei series, respectively.}\label{fig:bivst}\efi
\end{figure}

The sample cross-extremogram of any pair of series exhibits a
similar pattern of slow decay as seen in the univariate sample
extremograms (we do not include these figures). 
Figure~\ref{fig:bivst} shows the sample cross-extremograms  for the
filtered series. For example, in the first row
of graphs, $(X_t)$ is the filtered FTSE and
$(Y_t)$ are the filtered S\&P, DAX and Nikkei, respectively.  
There are signs of various types of cross-extremal dependence in the 
filtered series.  The spikes at lag zero (except between the Nikkei and S\&P)
indicate the strong extremal dependence of the multiplicative shocks, 
In the second row,
there is evidence of significant extremal dependence at
lag one for each sample cross-extremogram: given the
S\&P has an extreme left tail event in a shock at time $t$
there will be a corresponding large left tail shock in the FTSE, the DAX and the Nikkei
at time $t=1$.  Given the dominance of the US stock market, one might
expect a carry-over effect of the shocks 
on the other exchanges on the next day.  Since only marginal GARCH
models were fitted to the data, 
it may not seem all that surprising that the filtered series exhibit
serial dependence.  
We should note, however, that the dependence in the shocks does not appear to last beyond one time lag.

\section{A Fourier analysis of extreme events}
Classical \tsa\ studies the second order properties of
stationary processes  in the time and \fre y domains. The latter
approach refers to  spectral (or Fourier) analysis of the \ts .  
We mentioned in Section~\ref{subsec:extremogram} that the 
extremograms $\gamma_A=\gamma_{AA}$ and $\rho_A=\rho_{AA}$ for a set $A$ bounded away
from zero are covariance and correlations \fct s of some stationary
\seq . Then it is possible to study the corresponding spectral
properties of the extremogram. Research in this direction was started
in \cite{davis:mikosch:2009c} and continued in
\cite{mikosch:zhao:2012}. We recall some of the results.
Throughout we assume that $(X_t)$ is a $\bbr^d$-valued strictly
stationary \regvary\ \seq\ with index $\alpha>0$ and $A$ is 
a $\mu_1$-continuity set. In the comments following
Theorem~\ref{thm:1} we mentioned the latter property also implies that
the sets $A\times (\ov \bbr_0^d)^h\times A$ are $\mu_{h}$-continuity sets
for $h\ge 0$, so the limits $\rho_{A}(h)$ exist. 
\par
Assuming that
$\rho_{A}$ is square summable, we  consider the
corresponding spectral density corresponding to $\rho_{A}$
(see Brockwell and Davis \cite{brockwell:davis:1991}, Chapter 4):
\beao
f_A(\la)= 1+ 2\sum_{h=1}^\infty \rho_{A}(h)\,\ex^{-i\,h\,\la}\,,\quad \la\in [0,\pi]\,.  
\eeao
A natural estimator of the spectral density is obtained if we replace
the quantities $\rho_{A}(h)$ by the sample versions 
\beao
\wt
\rho_{A}(h)=\dfrac{\frac m n \sum_{t=1}^{n-h} (I_{a_m^{-1}X_t\in A}-p_0)
(I_{a_m^{-1}X_{t+h}\in A}-p_0)}{\wh P_m(A)}\,,\quad h\ge 1\,.
\eeao
where $p_0=P(a_m^{-1} X\in A)$:\footnote{The centering of the indicator \fct  s
with their expectation $p_0$ is crucial for deriving \asy\ theory.
In applications, $I_{nA}(\la)$ is typically evaluated at the Fourier
\fre ies $\w_j(n)=2\pi j/n\in (0,\pi)$ and since 
$\sum_{h=1}^n\ex^{-i\,h\,\w_j(n)}=0$, centering in $I_{nA}(\w_j(n))$
is not needed.}
\beao
\wt I_{nA}(\la) = 1+ 2\sum_{h=1}^\infty \wt \rho_{A}(h)\,\ex^{-i\,h\,\la}
= \dfrac{\dfrac mn\Big|\sum_{t=1}^n (I_{\{a_m^{-1} X_t\in
    A\}}-p_0)\,\ex^{i\,h\,\la}\Big|^2}
{\dfrac mn \sum_{t=1}^n I_{\{ a_m^{-1} X_t\in A\} }}=\dfrac{I_{nA}(\la)}{\wh P_m(A)}\,.
\eeao
We will refer to $I_{nA}$ and its standardized version $\wt I_{nA}$ as
{\em \per\ of the extreme event} $a_m \,A$. Indeed, if we replaced the
normalization $m/n$ by $1/n$,  $I_{nA}$ is the \per\ of the \seq\ of
centered indicators $(I_{\{a_m^{-1} X_t\in A\}}-p_0)$. These \seq s constitute a 
triangular array of row-wise strictly stationary \seq s for which
standard \asy\ theory for the \per\ does not apply; for an \asy\ 
theory of the \per\ of a strictly stationary linear process, see
Brockwell and Davis \cite{brockwell:davis:1991}, Chapter 10. However, 
the \per\ of extreme events shares some of the basic 
properties of the \per , as the following results from
Mikosch and Zhao \cite{mikosch:zhao:2012} show.

\bth\label{thm:2}
Let $(X_t)$ be an $\bbr^d$-valued strictly
stationary \regvary\ \seq\ with index $\alpha>0$ satisfying 
condition {\rm (M)}, $A\subset \ov\bbr_0^d$ be   
a $\mu_1$-continuity set and $\sum_{h\ge 1}\rho_A(h)<\infty$. 
\begin{itemize}
\item 
Assume $\la\in (0,\pi)$ is fixed and $\w_n=2\pi j_n/n$, $j_n\in \bbz$,
is any \seq\ of Fourier \fre ies \st\ $\w_n\to \la$. Then
\beao
\lim_{\nto} EI_{nA}(\la) =\lim_{\nto} EI_{nA}(\w_n)=\mu_1(A)\,f_A(\la)\,, 
\eeao
\item
Assume in addition that the \seq s  
$(m_n)$, $(r_n)$ from {\rm (M)} also satisfy  
the growth conditions $(n/m) \alpha_{r_n} \to 0$, and $m_n=o(n^{1/3})$. 
Let $(E_i)$ be a \seq\ of iid standard exponential \rv s.  
Consider any  fixed \fre ies  
$0<\lambda_1<\cdots< \lambda_N <\pi$ for some $N\ge 1$.  
Then the following relations hold:  
\beao  
\big(I_{nA}(\lambda_i)\big)_{i=1,\ldots,N}&\std&  
\mu_1(A)\,\big(f_A(\la_i) E_i\big)_{i=1,\ldots,N}\,,\quad \nto\,,\\  
\big(\wt I_{nA}(\lambda_i)\big)_{i=1,\ldots,N}&\std&  
\big(f_A(\la_i) E_i\big)_{i=1,\ldots,N}\,,\quad \nto\,.  
\eeao  
Consider any distinct Fourier \fre ies $\w_i(n)\to \la_i\in (0,\pi)$  
as $\nto$, $i=1,\ldots,N$. The limits $\la_i$ do not have to be distinct. Then  
the following relations hold:  
\beao  
\big(I_{nA}(\w_i(n))\big)_{i=1,\ldots,N}&\std&  
\mu_1(A)\,\big(f_A(\la_i) E_i\big)_{i=1,\ldots,N}\,,\quad \nto\,,\\  
\big(\wt I_{nA}(\w_i(n))\big)_{i=1,\ldots,N}&\std&  
\big(f_A(\la_i) E_i\big)_{i=1,\ldots,N}\,,\quad \nto\,.\\  
\eeao  
\end{itemize}  
\ethe
These properties are very similar to those of a strictly stationary
weakly dependent \seq .
The \asy\ independence of the \per\ of the extreme event $a_mA$
and consistency in the mean of $I_{nA}(\la)$ 
give raise to the hope that pointwise consistent smoothed 
\per\ estimation of the spectral
density $f_A(\la)$ is possible.
\par
\begin{figure}[htbp]  
\centerline{  
\includegraphics[height=7cm,width=10cm]{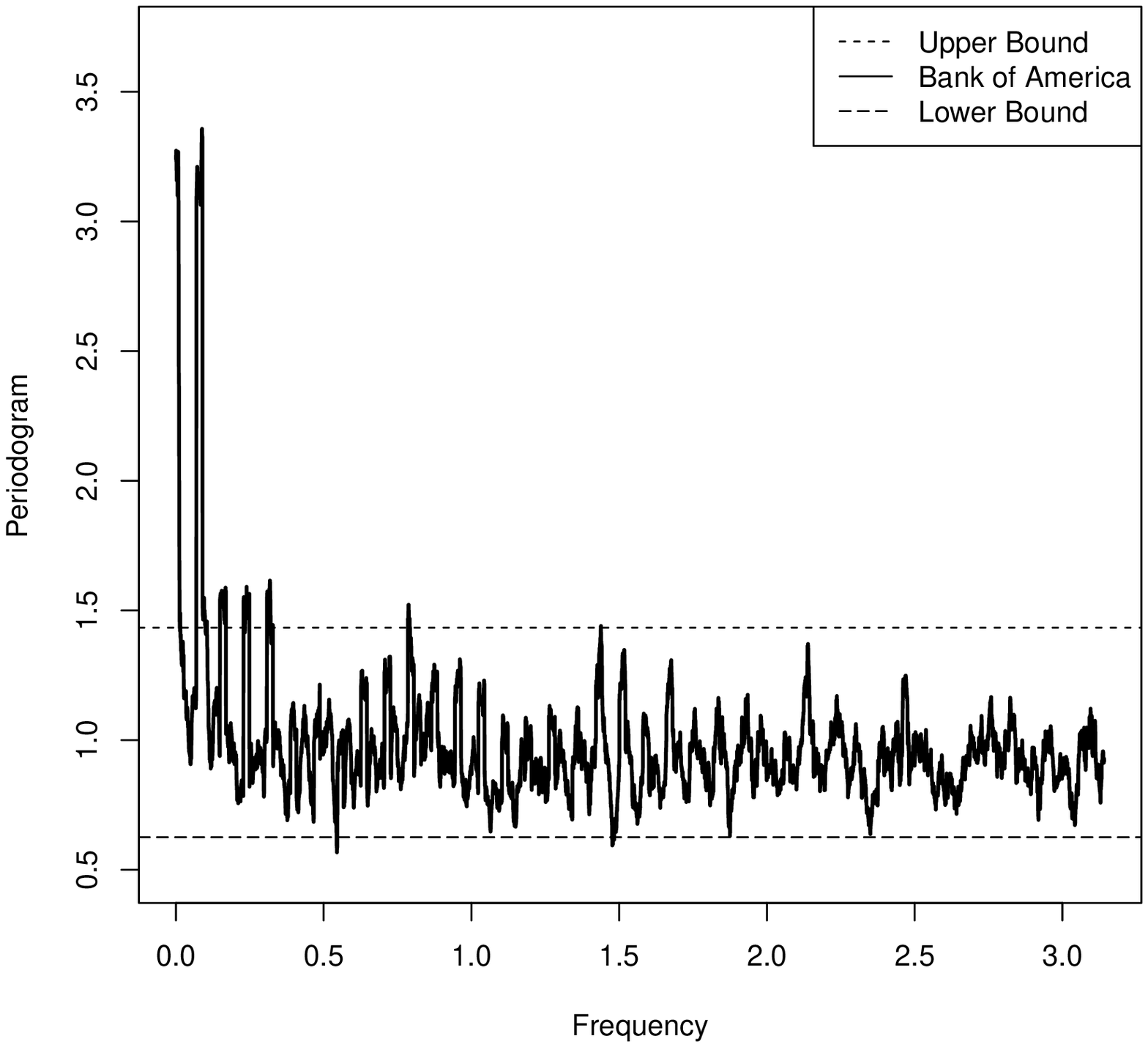}}  
\bfi{\em The smoothed periodogram (corresponding to the losses) for $31,757$ $5$-minute
  log-returns of Bank of America stock prices with 
Daniell window, $s_n=52$. The simultaneous 
confidence bands are constructed by taking   
the 97.5\% quantile of the maxima and the 2.5\% quantile of the minima 
over the Fourier \fre ies calculated from the smoothed \per s of 10 000 random permutations of the data. 
If the data were iid, permutations would not change the dependence structure. 
The fact that the \per\ is outside the confidence bands at various \fre ies 
indicates that there is significant extremal dependence in the data. 
The peaks at various \fre ies show that there are cycles of extremal behavior in the data. These cycles cannot be detected by autocorrelation plots of the data, their absolute values or squares.}
\efi
\label{fig:x}  
\end{figure}

For a fixed  \fre y $\la\in (0,\pi)$ define  
\beao  
\la_{0} = \min\{2\pi j/n : 2\pi j/n \geq  
\lambda\}\,,\quad\mbox{and}\quad\lambda_{j} = \lambda_{0} +2\pi  
j/n\,,\quad |j|\le s\,.  
\eeao  
(We suppress the dependence of $\la_j$ on $n$.)  
Assume that $s=s_n\to\infty$ and $s_n/n\to 0$  
as $\nto$. Consider the non-negative weight \fct\  
$(w_n(j)) _{|j|\le s}$  
satisfying the  
conditions  
\beao
& \sum_{|j|\leq s} w_n (j) =1 \quad\mbox{and}\quad  
\sum_{|j|\leq s} w_n^2 (j) \to 0 \mbox{ as } n \to \infty\,.  
\eeao   
Introduce the corresponding {\em smoothed \per }  
\beao  
\wh {f}_{nA} (\lambda) = \sum_{|j| \le s_n} w_n(j)\, I_{nA} (\lambda_{j}),  
\eeao  
Under the conditions of the second item in Theorem~\ref{thm:2} and
some further restrictions on the growth of $(m_n)$ and $(\alpha_h)$ 
the following limit relations hold for a fixed \fre y $\la\in
(0,\pi)$,  
\beao
\wh f_{nA}(\la)\stackrel{L^2}{\rightarrow} \mu_1(A)\,f_A(\la)
\quad\mbox{and} 
\quad  \wh f_{nA}(\la)/\wh P_m(A) \stp f_A(\la)\,.
\eeao  
In Figure~\ref{fig:x} we illustrate how the smoothed \per\ $\wh f_{n,(-\infty,-1)}$ 
works for  5-minute log-returns of Bank of America stock prices.
We choose a simple Daniell window with $w_n(j)=1/(2s_n+1)$ and $s_n=52$.

\end{document}